\documentclass[12pt,a4paper,reqno]{amsart}

\usepackage{graphicx}
\usepackage{footmisc}  
\usepackage{footnote}  
\makesavenoteenv{tabular}
\makesavenoteenv{table}
\usepackage{mathtools}
\usepackage{url} 
\usepackage[hidelinks,breaklinks, urlcolor=black]{hyperref} 

\usepackage[headings]{fullpage} 
\usepackage{time}  
\usepackage[english]{babel}  
\usepackage{enumerate}  
 
\usepackage{bm} 
\usepackage{mathrsfs}

\usepackage[full]{textcomp} 
\usepackage{newtxtext} % osf for text, not math
\usepackage{cabin} % sans serif 
\usepackage{zlmtt}% latin modern typewriter
\usepackage[bigdelims,vvarbb]{newtxmath} % bb from STIX
\usepackage[cal=boondoxo]{mathalfa} % mathcal

\usepackage{microtype} %   
\newcommand{\dx}{\mathrm{d}} 
\newcommand{\A}{\mathcal{A}} 
\newcommand{\B}{\mathcal{B}} 
\newcommand{\CC}{\mathcal{C}}

\newcommand{\C}{\mathbb{C}}
\newcommand{\Q}{\mathbb{Q}}
\newcommand{\N}{\mathbb{N}}

\newcommand{\Z}{\mathbb{Z}}
\newcommand{\G}{\mathfrak{G}}
\newcommand{\second}{{\prime\prime}}

\newcommand{\Odi}[1]{\Odip{}{#1}}

\newcommand{\Odip}[2]{\mathcal{O}_{#1}\left(#2\right)}

\newcommand{\odip}[2]{{o}_{#1}\left(#2\right)}
\newcommand{\odi}[1]{\odip{}{#1}}

\allowdisplaybreaks

\renewcommand{\qedsymbol}{$\square$}

\newtheoremstyle{sltheorems}% name
{10pt}%      Space above
{6pt}%      Space below
{\slshape}%         Body font
{}%         Indent amount (empty = no indent, \parindent = para indent)
{\bfseries}% Thm head font
{.}%        Punctuation after thm head
{.5em}%     Space after thm head: " " = normal interword space;
{\thmname{#1}\thmnumber{ #2}\thmnote{ (#3)}}
\theoremstyle{sltheorems}

\newtheoremstyle{remark}% name
{10pt}%      Space above
{6pt}%      Space below
{\rm} %         Body font
{}%         Indent amount (empty = no indent, \parindent = para indent)
{\bfseries}% Thm head font
{.}%        Punctuation after thm head
{.5em}%     Space after thm head: " " = normal interword space;
{\thmname{#1}\thmnumber{ #2}\thmnote{ (#3)}}
 \theoremstyle{remark}

\usepackage{caption} 
\allowdisplaybreaks
\makeatletter
\def\subsubsection{\@startsection{subsubsection}{3}%
  \z@{.3\linespacing\@plus.5\linespacing}{-.5em}%
  {\normalfont\bfseries}}
\makeatother

\begin{document} 

\title[Computation of the Euler-Kronecker  constants]{Efficient computation of the Euler-Kronecker \\constants of 
prime cyclotomic fields}
\author{Alessandro Languasco}
\date{} 

\begin{abstract}  
We introduce a new algorithm, which is faster and requires 
less computing resources than the ones previously known,  to compute the
Euler-Kronecker constants  $\G_q$ for the prime cyclotomic fields $ \Q(\zeta_q)$,
where $q$ is an odd prime and $\zeta_q$ is a primitive $q$-root of unity.
With such a new algorithm we evaluated $\G_q$ and $\G_q^+$, 
where $\G_q^+$ is the Euler-Kronecker constant of the maximal real subfield of 
$\Q(\zeta_q)$,  for some very large primes $q$ thus obtaining two new negative values of $\G_q$:
 $\G_{9109334831}=  -0.248739\dotsc$  and $\G_{9854964401}=  -0.096465\dotsc$ 
We also evaluated  $\G_q$ and $\G^+_q$  for every odd prime
$q\le 10^6$, thus enlarging the size of the previously known range for $\G_q$ and $\G^+_q$.
Our method also reveals that difference $\G_q - \G^+_q$ can be computed in a much simpler 
way than both its summands, see Section \ref{EKdiff-formula}.
Moreover, as a by-product, we  also computed 
$M_q=\max_{\chi\ne \chi_0} \vert L^\prime/L(1,\chi) \vert $  
 for every odd prime
$q\le 10^6$, where   $L(s,\chi)$ are the Dirichlet $L$-functions, $\chi$ run over the 
non trivial Dirichlet characters mod $q$ 
and $\chi_0$ is the  trivial Dirichlet character mod $q$.
 As another by-product of our computations, we will also 
 provide more data on the generalised Euler constants
 in arithmetic progressions.    
\end{abstract} 
\maketitle 
\section{Introduction} %
  Let $K$ be a number field and let $\zeta_K(s)$ be its Dedekind zeta-function. 
 It is a well known fact that $\zeta_K(s)$ has a simple pole at $s=1$; 
  writing the expansion of $\zeta_K(s)$ near $s=1$  as
  \[
  \zeta_K (s)  = \frac{c_{-1}}{s-1} + c_0 + \Odi {s-1},
  \]
 the \emph{Euler-Kronecker constant of} $K$ is defined as    
  \[
\lim_{s\to 1} \Bigl(\frac{\zeta_K (s)}{c_1} - \frac{1}{s-1}\Bigr)= \frac{c_0}{c_{-1}}.
  \]
  In the special case  in which  $K=\Q(\zeta_q)$ is a prime   cyclotomic  field, 
  where $q$ is an odd prime and $\zeta_q$ is a primitive $q$-root of unity,
we have that the Dedekind zeta-function  satisfies
  $\zeta_{\Q(\zeta_q)} (s)= \zeta(s) \prod_{\chi \neq \chi_0} L(s,\chi) $,
  where $\zeta(s)$ is the Riemann zeta-function, 
$L(s,\chi)$ are the Dirichlet $L$-functions, $\chi$ runs over the non trivial Dirichlet 
characters mod $q$ 
and $\chi_0$ is the  trivial Dirichlet character mod $q$.
By logarithmic differentiation, we immediately get  that
the \emph{Euler-Kronecker constant
for the prime cyclotomic field $\Q(\zeta_q)$} is
\begin{equation}
\label{EKq-def}
\G_{q}
: =
\gamma
+
\sum_{\chi \neq \chi_0} \frac{L^\prime}{L}(1,\chi),
\end{equation}
where
$\gamma$ is
the Euler-Mascheroni constant.
Sometimes the quantity $\G_q$ is denoted as  $\gamma_q$
but this conflicts with notations used in literature.
 
An extensive study of the properties of $\G_{q}$  was recently started by 
Ihara \cite{Ihara2006,Ihara2008} and continued by many others; 
here we are mainly interested in computational problems
involving $\G_q$ and hence we just recall the paper by Ford-Luca-Moree \cite{FordLM2014}.
We introduce a  new method to compute the
 Euler-Kronecker constants of prime cyclotomic fields
 which is faster and uses less computing resources than 
 the ones previously known.
The new algorithm requires
 the values of the \emph{generalised gamma function} $\Gamma_1$
 at some rational arguments  $a/q\in(0,1)$. Such a function  has, for $q$ large and $a=\odi{q}$, 
 an order of magnitude exponentially smaller than the ones previously used
 to determine $\G_q$, see Section \ref{EKq-comput} below. 
 Moreover, the Fast Fourier Transform (FFT) used in this new approach
 allows a \emph{decimation in frequency strategy}\footnote{We use here this nomenclature 
 since it is standard in the literature on the Fast Fourier Transform, but it could be 
 translated in number theoretic language using suitable properties of Dirichlet characters.}  
 that leads to  gaining a factor $1/2$ 
 in the quantity of needed precomputation operations, in the length of the involved transforms
  and in their memory usage. Our algorithm uses the formulae in Section \ref{dist-parity}  below
  that, according to  Deninger \cite{Deninger1984} and  Kanemitsu \cite{Kanemitsu1989}, 
  were first proved in 1883 by Berger \cite{Berger1883}
  and in 1929 by Gut \cite{Gut1929}.

Another interesting quantity related to $\G_q$ is
the Euler-Kronecker constant $\G_q^+$ for $\Q(\zeta_q+\zeta_q^{-1})$, the maximal real subfield
of $\Q(\zeta_q)$. According to eq.~(10) of Moree \cite{Moree2018}  it is defined as
\begin{equation}
\label{EKq+-def}
\G_q^+ : = \gamma 
+
\sum_{\substack{\chi \neq \chi_0\\ \chi\, \textrm{even}}} 
\frac{L^\prime}{L}(1,\chi).
\end{equation}
In Section \ref{EKdiff-formula} we will give a formula that leads a to  a direct evaluation of $\G_q^+$
in terms of some special functions values attained at some rationals  $a/q\in(0,1)$.
Moreover, in Section \ref{EKdiff-formula} we will use  the previously proved relations
to see why the quantity $\G_q-\G_q^+$   is much easier to compute than  both its summands.

During such computations,  as a by-product, we also evaluated the related quantity
\begin{equation}
\label{Mq-def}
M_q:=\max_{\chi \neq \chi_0} \Bigl\vert\frac{L^\prime}{L}(1,\chi) \Bigr\vert,
\end{equation}
see Section \ref{Mq-eval}.
Other quantities related to $\G_q$ are
  the \emph{generalised Euler constants in arithmetic progressions},
  sometimes also called \emph{Stieltjes constants in arithmetic progressions},
  denoted as  $\gamma_k(a,q)$, $k\in \N$,  $q \ge 1$,
 $1\leq a \leq q$, which are defined by
  \begin{align}
\notag
 \gamma_k(a,q) &:= 
 \lim_{N\to+\infty}
\Bigl(
\sum_{\substack{0<m\le N  \\ m  \equiv a \bmod q}} \frac{(\log m)^k}{m}
- 
\frac{(\log N)^{k+1}}{q(k+1)}
\Bigr)
\\ &
  \label{gammakaq-def}
 =
 - \frac{1}{q}
 \Bigl(
 \frac{(\log q)^{k+1}}{k+1} +
 \sum_{n=0}^{k} 
 \binom{k}{n} 
 (\log q)^{k-n} \psi_n\bigl(\frac{a}{q}\bigr)
 \Bigr),
 \end{align}
see eq.~(1.3)-(1.4) and (7.3) of  Dilcher \cite{Dilcher1992},
  where 
  \begin{equation}
  \label{gen-digamma-def}
   \psi_n(z) := -\gamma_n - \frac{(\log z)^n}{z} 
   - \sum_{m=1}^{+\infty}
   \Bigl( 
   \frac{(\log (m+z))^n}{m+z} -    \frac{(\log m)^n}{m}
   \Bigr) 
  \end{equation}
for $n\in \N$ and $z\in \C\setminus\{0,-1,-2,\dotsc \}$, $\psi_n(1) = -\gamma_n$, 
and
the \emph{generalised Euler constants} $\gamma_n$ are defined as
\begin{equation}
\label{gamman-def}
\gamma_n :=
\lim_{N\to+\infty}
\Bigl(
\sum_{j=1}^{N}
\frac{(\log j)^n}{j} - \frac{(\log N)^{n+1}}{n+1}
\Bigr)
=
 \sum_{m=1}^{+\infty}
   \Bigl( 
   \frac{(\log m)^n}{m} -    \frac{(\log (m+1))^{n+1}-(\log m)^{n+1}}{n+1}
   \Bigr),
\end{equation}
by, \emph{e.g.}, eq.~(3)-(4) of Bohman-Fr\"oberg \cite{BohmanF1988}.
Remark that $\gamma_0=\gamma$.
It is worth recalling that the functions $\psi_n(z)$ occur in Ramanujan's second
notebook, see \cite[Chapter 8, Entry 22]{Berndt1985}.

The quantities in \eqref{gammakaq-def} and, as we will see in sections 
\ref{general-case}-\ref{dist-parity}
below, the one  in \eqref{EKq-def},  are hence connected with 
the values of $\psi_n$, $n\ge1$, which is the logarithmic derivative of $\Gamma_n$, a generalised 
Gamma function, see Deninger \cite{Deninger1984}, Dilcher \cite{Dilcher1994} and 
Katayama \cite{Katayama2013},
whose definition for $n=1$ is given in Section \ref{even-case-Deninger}. 
In some sense we can say that  the $\psi_n$-functions, $n\ge 1$, are the analogue of the 
usual \emph{digamma} function. In the following we will denote as $\psi$ the standard 
digamma function $\Gamma^\prime/\Gamma$; we also remark that
it can be represented as the function $\psi_0$ defined in \eqref{gen-digamma-def}. 
 
The paper is organised as follows. In Sections \ref{general-case}-\ref{dist-parity} 
we will give the  derivations  of the main formulae we will need in the computations;
such proofs are classical and are based on the  functional equation for 
the Dirichlet $L$-functions, see Cohen's books  \cite{Cohen2007a}-\cite{Cohen2007}, for instance. 
Other useful references  for this part are  the
papers of Deninger \cite{Deninger1984} and Dilcher \cite{Dilcher1992,Dilcher1994}.

In Section \ref{EKq-comput} we will see how to implement the formulae of the previous
two sections, starting form a straightforward application of the definitions 
\eqref{EKq-def}-\eqref{EKq+-def} of $\G_q$
and $\G_q^+$;
then we will  compare the Ford-Luca-Moree approach, based on 
 the formulae of Section \ref{general-case},
with our new procedure, based on  the formulae of Section \ref{dist-parity}.
In particular 
we will see, in both cases, how to insert  the Fast Fourier Transform  and we will discuss
their precisions, computational costs and memory usages. 
In Section \ref{Mq-eval} we describe how to compute $M_q$.  
Section \ref{gen-euler-arit-prog}  is devoted to  provide more data on the generalised 
Euler constants
 in arithmetic progressions.   
 Finally, in Section \ref{tables} we will provide several tables containing
a comparison scheme of the different implementations and the computational results
and running times of the different approaches. At the bottom of the paper we inserted two 
colored
scatter plots for the  normalised values
of $\G_q$ and $\G_q^+$ for every prime $q$, $3\le q \le 10^6$
and two scatter plots about  $M_q$ and its normalised values for the same set 
of primes. %

We finally remark that some of the ideas presented here will also be used in
a joint work with Pieter Moree, Sumaia Saad Eddin and Alisa Sedunova
on the computation of the Kummer ratio of the class number for prime cyclotomic fields,
see \cite{LanguascoMSS2019}.

\medskip 
\textbf{Acknowledgements}. 
Some of the calculations here described were performed
using the University of Padova Strategic Research Infrastructure Grant 2017:
``CAPRI: Calcolo ad Alte Prestazioni per la Ricerca e l'Innovazione'',
 \url{http://capri.dei.unipd.it}.
I also wish to thank Karim Belabas and Bill Allombert (University of Bordeaux) for a couple of key suggestions
about  {\tt libpari} and {\tt gp2c} and  Luca Righi (University of Padova) for his help in developing the quadruple precision
versions of the fft-programs, in designing the parallelised precomputations  and in 
organising the  use of the cluster of the 
Dipartimento di Matematica  ``Tullio Levi-Civita'', 
\url{http://computing.math.unipd.it/highpc}, and the use of CAPRI.

\section{The Ford-Luca-Moree method}
\label{general-case} 
Recall that $q$ is an odd prime.
If we do not restrict to  Dirichlet characters of a prescribed parity, we can    
use eq.~(6.1) and (7.4) of Dilcher \cite{Dilcher1992}, as in Ford-Luca-Moree, see eq.~(3.2) 
in \cite{FordLM2014}. 
In fact eq.~(6.1) of \cite{Dilcher1992} gives
\[
L^\prime(1,\chi)
=
- 
\sum_{a=1}^{q-1}  \chi(a)\gamma_{1}(a,q),
\]
where $\gamma_{1}(a,q)$ is defined in \eqref{gammakaq-def} which, for $k=1$, becomes
\begin{equation*} 
\gamma_{1}(a,q) 
=
- \frac{1}{q}
\Bigl(
\frac{1}{2}(\log q)^{2}
+ \log q\  
\psi\bigl(\frac{a}{q}\bigr) 
+ \psi_{1}\bigl(\frac{a}{q}\bigr)
\Bigr),
\end{equation*}
 for any $q\geq 1$ and $1\leq a \leq q$,
where 
$\psi,\psi_1$ are defined in \eqref{gen-digamma-def}.
Again using \eqref{gen-digamma-def}, we  define 
\begin{equation}
\label{T-def}
T(x) 
:= 
 \gamma_1 + \psi_{1}(x)
 =
 -\frac{\log x}{x}
-
\sum_{m=1}^{+\infty}
\Bigl(
\frac{\log(x+m)}{x+m} - \frac{\log m}{m}
\Bigr),
\end{equation}
and, specialising \eqref{gamman-def}, we also have 
\begin{equation} 
\notag
\gamma_{1}  =
\lim_{N\to+\infty}
\Bigl(
\sum_{j=1}^{N}
\frac{\log j}{j} - \frac{(\log N)^2}{2}
\Bigr) 
=-0.0728158454835\dotsc 
\end{equation}
To compute $\gamma_1$  and   similar constants with a very large precision, see 
Section \ref{gammak-aq-comp} below.
We also remark here that the rate of convergence of the series in \eqref{T-def}
is, roughly speaking, about $(\log m) /m^2$.
Recalling  now eq.~(3.1) of \cite{FordLM2014}, \emph{i.e.},
\begin{equation}
\label{general-2}
L(1,\chi)
=
-
\frac{1}{q}
\sum_{a=1}^{q-1} \chi(a)\ 
\psi \bigl(\frac{a}{q}\bigr), 
\end{equation}
by the orthogonality of Dirichlet characters and \eqref{general-2}, we
obtain eq.~(3.2) of \cite{FordLM2014}, \emph{i.e.},
\begin{equation*}
L^\prime(1,\chi) 
=
- (\log q) L(1,\chi)
+
\frac{1}{q}
\sum_{a=1}^{q-1}
\chi(a)\ 
T\bigl(\frac{a}{q}\bigr),
\end{equation*}
where $T(x)$ is defined in \eqref{T-def}
(pay attention to the change of sign in \eqref{T-def} with respect to eq.~(3.2) of \cite{FordLM2014}).
Summarising, we finally get
\begin{equation}
\label{general-simplified}
\sum_{\chi \neq \chi_0}\frac{L^\prime}{L}(1,\chi) 
=
-(q-2)\log q
-
\sum_{\chi \neq \chi_0}
\frac{\sum_{a=1}^{q-1} \chi(a) \ T(a/q)}
{\sum_{a=1}^{q-1} \chi(a)\ 
\psi (a/q) 
}.
\end{equation}

Formula \eqref{general-simplified} is the one used in the paper 
by Ford-Luca-Moree \cite{FordLM2014}. We will now explain how we can compute $\G_q$ via \eqref{EKq-def}   using the 
values of the two special functions $\psi$ and $T$,  together with the values of the non trivial 
Dirichlet characters mod $q$.  

From a computational point of view it is clear that in \eqref{general-simplified} 
we first have to evaluate  $T(a/q)$ and 
 $ \psi(a/q)$ for every $1\leq a \leq q-1$.
For the $\psi$-values
we can rely  on the PARI/GP function {\tt psi} or, if less precision is sufficient,
we can use the analogous function included in  GSL, the \emph{gnu scientific library} \cite{GSL2018}\footnote{GSL 
provides just a \emph{double precision}  (in the sense of the C programming language precision) 
version of $\psi$; hence this is faster, but less accurate, 
than the computation of the $\log \Gamma$-values needed in the procedure described in the next section. If 
we use  PARI/GP to precompute and store the $\psi$-values, then the costs of the precomputation  and 
the input/output part of the FFT step have to be doubled, see Table \ref{table1}.}. 
For  computing
the $T$-values, a task for which there are no pre-defined functions 
in any software libraries we know,  we can use the summing function
{\tt sumnum} of PARI/GP; this is the most time-consuming step of the procedure.
Using the FFT algorithm 
to perform the sums over $a$,
it is easy to see that computing $\G_q$ via \eqref{general-simplified}
has a computational cost of $\Odi{q\log q}$ arithmetical operations together with the cost of computing 
$q-1$ values of the $\psi$ and $T$ functions.
For more details see Section \ref{EKq-comput}.

\section{Another method: distinguishing Dirichlet characters' parities}
\label{dist-parity}
\subsection{Primitive odd Dirichlet character case.}
\label{new-odd}
Recall that $q$ is an odd prime,
let $\chi\neq \chi_0$ be a primitive odd Dirichlet character mod $q$ and
let $\tau(\chi):= \sum_{a=1}^q \chi(a)\, e(a/q)$, $e(x):=\exp(2\pi i x)$, be the Gau\ss\ sum 
associated with $\chi$.
The functional equation  for $L(s,\chi)$, see, \emph{e.g.},  the proof of Theorem  3.5 of 
 Gun-Murty-Rath \cite{GunMR2011}, gives
\[
L(s,\chi)
=
\frac{1}{\pi i} 
\Bigl( \frac{2 \pi}{q} \Bigr)^s
\Gamma(1-s)
\frac{\tau(\chi)}{\sqrt{q}}
\cos\Bigl( \frac{\pi s}{2} \Bigr)
L(1-s,\overline{\chi})
\]
and hence
\[
\frac{L^\prime}{L}(s,\chi)
=
\log  \Bigl( \frac{2 \pi}{q} \Bigr)
- \frac{\Gamma^\prime}{\Gamma}(1-s)
-  \frac{\pi }{2} \tan\Bigl( \frac{\pi s}{2} \Bigr)
- \frac{L^\prime}{L}(1-s,\overline{\chi}),
\]
which, evaluated  at $s=0$, gives
\begin{equation}
\label{derlog-01-rel}
\frac{L^\prime}{L}(0,\chi)
=
\log  \Bigl( \frac{2 \pi}{q} \Bigr)
+\gamma
- \frac{L^\prime}{L}(1,\overline{\chi}).
\end{equation}
By the Lerch identity about values of the Hurwitz zeta-function, 
see, \emph{e.g.},
Proposition 10.3.5 of Cohen \cite{Cohen2007},
and the orthogonality of Dirichlet characters, we get
\begin{align}
\notag
L^\prime(0,\chi)
&=
-
\log q \sum_{a=1}^{q-1} \chi(a)\Bigl(\frac{1}{2}-\frac{a}{q}\Bigr)
+
\sum_{a=1}^{q-1}  \chi(a) \log\Bigl(\Gamma\bigl(\frac{a}{q}\bigr)\Bigr)
\\
\notag
&
=
\frac{\log q}{q}  \sum_{a=1}^{q-1} a \chi(a) 
+
\sum_{a=1}^{q-1}  \chi(a) \log\Bigl(\Gamma\bigl(\frac{a}{q}\bigr)\Bigr)
\\
\label{Lprime0}
&
=
- (\log q) L(0,\chi)
+
\sum_{a=1}^{q-1}  \chi(a) \log\Bigl(\Gamma\bigl(\frac{a}{q}\bigr)\Bigr),
\end{align}
since, see  Proposition 9.5.12 and Corollary 10.3.2  of Cohen \cite{Cohen2007}, we have
\begin{equation}
\label{L0}
L(0,\chi)
=
- B_{1,\chi}
: = 
 - \frac{1}{q} \sum_{a=1}^{q-1}  a \chi(a) 
 \ne 0,
\end{equation}
where $B_{1,\chi}$ is  the first $\chi$-Bernoulli number which is non-zero since $\chi$ is odd. 
 
Summarising, by \eqref{derlog-01-rel}-\eqref{L0}, we obtain 
\begin{align} 
\sum_{\chi\, \textrm{odd}} 
\frac{L^\prime}{L}(1,\chi) 
\label{odd-2}
&=
\frac{q-1}{2}
\bigl(
\gamma 
+
\log ( 2 \pi )
\bigr)
+ 
\sum_{\chi\, \textrm{odd}} 
\frac{1}{B_{1,\overline{\chi}} }
 \sum_{a=1}^{q-1}   \overline{\chi}(a)\log\Bigl(\Gamma\bigl(\frac{a}{q}\bigr)\Bigr).
\end{align}

From a computational point of view,
in \eqref{odd-2} we need to compute the $\log\Gamma$-values; to do so 
 we can rely  on an internal  PARI/GP function
 or, if less precision is sufficient,
we can use the analogous function included in the C programming language.
We  remark that,
for $x\to 0^+$,  $\log\bigl(\Gamma(x)\bigr)\sim \log(1/x)$ and
$\psi(x)\sim  -1/x $; hence for $q$ large and $a=\odi{q}$,  the values of 
$\log\bigl(\Gamma(a/q)\bigr)$ are exponentially smaller than the ones of
$\psi(a/q)$. 
Moreover, to compute the first $\chi$-Bernoulli number $B_{1,\overline{\chi}}$, defined in \eqref{L0},
we just need an integral sequence.   

\subsection{Primitive even Dirichlet character case.}
\label{even-case-Deninger}
Recall that $q$ is an odd prime.
Assume now that  $\chi\neq \chi_0$ is a primitive even Dirichlet character mod $q$.
We follow Deninger's notation in \cite{Deninger1984} by writing 
$R(x)= - \frac{\partial^2}{\partial s^2} \zeta(s,x) \vert _{s=0}
=\log (\Gamma_1(x))$, $x>0$, where $\zeta(s,x)$ is the Hurwitz zeta function, 
$s\in \C\setminus\{1\}$.
By eq.~(3.5)-(3.6) of  \cite{Deninger1984} we have
\begin{equation}
\label{even-1}
L^\prime(1,\chi)
=
(\gamma + \log(2\pi)) L(1,\chi)
+
\frac{\tau(\chi)}{q}
\sum_{a=1}^{q-1} \overline{\chi}(a)\ 
R\bigl(\frac{a}{q}\bigr),
\end{equation}
where, see eq.~(2.3.2) of \cite{Deninger1984}, 
the $R$-function can be expressed for every $x>0$ by
\begin{align}
\label{R-def}
R(x) &:=  -\zeta^{\prime\prime}(0) - S(x),\\
\label{S-def}
S(x) 
&:=
2 \gamma_1 x +(\log x)^{2} 
+ 
    \sum_{m=1}^{+\infty}
\Bigl(
\bigl(\log (x+m)\bigr)^{2} - (\log m)^{2} -2x \frac{\log m}{m}
\Bigr) .
\end{align}
It is worth recalling that 
comparing \eqref{R-def}-\eqref{S-def}  with \eqref{T-def}, we see that 
$\psi_{1}(x)=R^\prime(x)/2$ (note the different definition of $\gamma_{1}$ on page 174 of  
Deninger's paper).
Using \eqref{gamman-def} we have $S(1)=0$ and  $R(1) =  -\zeta^{\prime\prime}(0)$.
An alternative definition of $S(x)$ for $x>0$, which will be useful during the computations,
is implicitly contained in eq.~(2.12) of Deninger \cite{Deninger1984}:
\begin{equation}
\label{S-alt-def}
S(x) =
2 \int_{0}^{+\infty} 
\Bigl( (x-1) e^{-t} + \frac{e^{-xt} - e^{-t}}{1-e^{-t}} \Bigr)\frac{\gamma + \log t}{t} \ \dx t,
\quad x>0.
\end{equation}

By the orthogonality of the Dirichlet characters, we  immediately get
\begin{equation}
\label{R-S-summed}
 \sum_{a=1}^{q-1} \overline{\chi}(a)\,
R (a/q )
=
-
\sum_{a=1}^{q-1} \overline{\chi}(a)\,
S (a/q).
\end{equation}
For $L(1,\chi)$, we use  formula (2) of Proposition 10.3.5 of Cohen
\cite{Cohen2007}
and the parity of $\chi$ to get
\begin{equation}
\label{even-2} 
L(1,\chi)
= 
2\frac{\tau(\chi)}{q}
\sum_{a=1}^{q-1} \overline{\chi}(a)\log\Bigl(\Gamma\bigl(\frac{a}{q}\bigr)\Bigr) , 
\end{equation}
since $W(\chi)=\tau(\chi)/q^{1/2}$ for even Dirichlet characters, 
see Definition 2.2.25 of Cohen \cite{Cohen2007a}.
Summarising,  using \eqref{even-1} and \eqref{R-S-summed}-\eqref{even-2},
if $\chi$ is an even Dirichlet character mod $q$, we finally get
\begin{align} 
\sum_{\substack{\chi \neq \chi_0\\ \chi\, \textrm{even}}} 
\frac{L^\prime}{L}(1,\chi) 
\label{even-simplified}
&= 
\frac{q-3}{2}
\bigl(
\gamma + \log(2\pi)
\bigr)
-\frac{1}{2}
\sum_{\substack{\chi \neq \chi_0\\ \chi\, \textrm{even}}} 
\frac{\sum_{a=1}^{q-1} \overline{\chi}(a) \ S(a/q)}
{\sum_{a=1}^{q-1} \overline{\chi}(a)\log\bigl(\Gamma(a/q)\bigr)}.
\end{align}

We remark that in \eqref{even-simplified}   we can reuse the $\log\Gamma$-values  
already needed in \eqref{odd-2}. For computing the $S$-values, a task for which there 
are no pre-defined functions 
in any software libraries we know, we can use the PARI/GP  functions
{\tt sumnum} and {\tt intnum}; this is the most time-consuming step of the procedure. 
We also remark that,
for $x\to 0^+$,  $S(x) \sim   (\log x)^2$ and
$T(x)\sim \log(1/x)/x$; hence for $q$ large and $a=\odi{q}$,  the  values of 
$S(a/q)$ are exponentially smaller than the corresponding ones of
$T(a/q)$.
 
Finally we  remark that for $S$ and $B_{1,\overline{\chi}}$ we just need  the summation
over half of the Dirichlet characters involved; hence  in both cases in their computation 
using the Fast Fourier Transform
we can implement the so-called \emph{decimation in frequency} strategy that allow us 
to improve
on both the speed and the memory usage of the actual computation, 
see Section \ref{FFT-technique} below.
Using the FFT algorithm to perform the sums over $a$,
it is easy to see that computing $\G_q$ via \eqref{odd-2} and \eqref{even-simplified}
has a computational cost of $\Odi{q\log q}$ arithmetical operations together with the cost of computing 
$q-1$ values of the $\log \Gamma$-function and $(q-1)/2$  decimated in frequency 
values of the $S$-function\footnote{An explanation for this fact can be found at the end of 
Section \ref{FFT-technique}.}.

\subsection{On $\G_q^+$: the constant attached to the maximal real subfield
of $\Q(\zeta_q)$.} 
It is a consequence of the computations   in this section  that 
the Euler-Kronecker constant $\G_q^+$ for $\Q(\zeta_q+\zeta_q^{-1})$, the maximal 
real subfield
of $\Q(\zeta_q)$, is directly connected with the $S$-function
since, by \eqref{EKq+-def}
and \eqref{even-simplified},
we have
\begin{equation}
\label{EK+-formula}
\G_q^+   
=
\frac{q-1}{2}\gamma + \frac{q-3}{2}\log(2\pi)
-
\frac{1}{2}
\sum_{\substack{\chi \neq \chi_0\\ \chi\, \textrm{even}}} 
\frac{\sum_{a=1}^{q-1} \overline{\chi}(a) \ S(a/q)}
{\sum_{a=1}^{q-1} \overline{\chi}(a)\log\bigl(\Gamma(a/q)\bigr)}.
\end{equation}
Hence  in this case  the relevant information 
is encoded in the $S$ and $\log \Gamma$ functions. 
Clearly  $\G_q^+$  can be obtained during the $\G_q$-computation
since it requires a subset of  the data needed for getting $\G_q$.
In Figure \ref{fig2} you can find a 
colored
scatter plot of its values for every $q$
prime, $3\le q \le 10^6$. 

Moreover, a direct evaluation of $\G_q^+$ via \eqref{EK+-formula}
allow us to use a \emph{decimation in frequency} strategy in the 
application of the FFT technique to evaluate the sums over $a$, see Sections
\ref{FFT-technique}-\ref{DIF-odd}.

\subsection{Regarding $\G_q-\G_q^+$.}
\label{EKdiff-formula}
By \eqref{EKq-def}-\eqref{EKq+-def},  \eqref{odd-2} and \eqref{EK+-formula} it is 
trivial to get that
\begin{equation}
\label{diff-equation}
\G_q-\G_q^+ 
= 
\sum_{\chi\, \textrm{odd}} \frac{L^\prime}{L}(1,\chi)
=
\frac{q-1}{2}\bigl(\gamma + \log(2\pi)\bigr)
+   
\sum_{\chi\, \textrm{odd}} 
\frac{1}{B_{1,\overline{\chi}} }
 \sum_{a=1}^{q-1}   \overline{\chi}(a)\log\Bigl(\Gamma\bigl(\frac{a}{q}\bigr)\Bigr).
\end{equation}

This reveals that, from a practical point of view, $\G_q-\G_q^+$ is  much easier to 
compute with respect to both $\G_q$ and $\G_q^+$: this 
not just because, as for $\G_q^+$, it requires a subset of  the data needed  for $\G_q$ 
but also because
it involves   just one special function, $\log \Gamma$, which is  directly 
available in many software libraries and in the C  programming language.

In this case too, a direct evaluation of $\G_q - \G_q^+$ via \eqref{diff-equation}
allow us to use a \emph{decimation in frequency} strategy in the 
application of the FFT technique to evaluate the sums over $a$, see Sections
\ref{FFT-technique}-\ref{DIF-odd}. Some computational data about this quantity
are also included in \cite{LanguascoMSS2019}.

\section{Comparison of methods, results and running times} 
\label{EKq-comput}

First of all we notice that PARI/GP, v.~2.11.4, has the ability to  generate
 the Dirichlet $L$-functions (and  many other $L$-functions)
and hence the computation of $\G_q$, $\G_q^+$ and $M_q$   can be  performed using
 \eqref{EKq-def}-\eqref{Mq-def} 
with few instructions of the gp scripting language.
This computation has a linear cost in the number of calls
of the {\tt lfun} function of PARI/GP and,  at least  
on our Dell Optiplex desktop machine,  it is slower than 
both the procedures we are about to describe.
 
Comparing \eqref{odd-2} and \eqref{even-simplified} with \eqref{general-simplified},
we see that in both cases we can rely on pre-defined functions 
to compute either the $\log\bigl(\Gamma(a/q)\bigr)$-values or  
the $\psi(a/q)$-values, $1\leq a \leq q-1$, and finally we have  to evaluate the $T$ and 
$S$ functions
respectively involved. We recall that, when taking $q$ very large, 
it is  relevant to know their order of magnitude for $x\to 0^+$; it is easy to verify that
$\log\bigl(\Gamma(x)\bigr)\sim \log(1/x)$,  $S(x) \sim   (\log x)^2$,  
$\psi(x)\sim  -1/x $ and $T(x)\sim \log(1/x)/x$. Hence 
for $x\to 0^+$,  we have that   $\log\bigl(\Gamma(x)\bigr)$ and  $S(x)$ 
are exponentially smaller
than $\psi(x)$ and $T(x)$; a fact that will lead to a more accurate result  
when using a fixed precision in the final step of the computation.
Another difference is  that,
for the odd Dirichlet characters, the first $\chi$-Bernoulli number in eq.~\eqref{odd-2}
does not involve  any special function, but  just an integral sequence. So
it seems reasonable to compare the following two approaches:
\begin{enumerate}[a)]
\item 
\label{T-approach}
use the $T$-series formulae and the $\psi$-values 
as in \cite{FordLM2014}; in this case we have two possible
alternatives to evaluate the $\psi$-function: using GSL (gaining
in speed but losing in precision) or using PARI/GP (with a
much better precision, but doubling the needed hard disk storage
and the number of input/output operations on the hard disk);
\item 
\label{S-approach}
use the $S$-function formulae for the even Dirichlet characters case and the first $\chi$-Bernoulli
number for the odd one; remark that
in both cases we have to evaluate a sum of  the $\log\Gamma$-values.
\end{enumerate}

 This way we can extend the computation performed
in \cite{FordLM2014}, not only because we are 
developing a different implementation of the same formulae, 
but also because we can solve the problem in an alternative way
which is faster, needs less computing resources, and uses functions having a much smaller 
order of magnitude,
see Table \ref{table1} for a summary of these facts. 
In the computation we will use the PARI/GP
scripting language to exploit its ability to accurately evaluate
the series and integrals involved in the definition of the $T$ and $S$ functions,
defined respectively in \eqref{T-def} and \eqref{S-def}-\eqref{S-alt-def}, via the  functions
{\tt sumnum} or {\tt intnum}.

  \subsection{Using the FFT algorithm}
  \label{FFT-technique}
 We also remark that the procedures \ref{T-approach})-\ref{S-approach}) trivially require a
quadratic number of arithmetical operations to perform  the computations in
\eqref{general-simplified}, \eqref{odd-2} and \eqref{even-simplified},
but this can be improved by using the  FFT algorithm and the following argument.
Focusing on  \eqref{general-simplified}, \eqref{odd-2} and \eqref{even-simplified}, we remark that,
since $q$ is prime, it is enough to get $g$, a primitive root of $q$,
and $\chi_1$, the Dirichlet character mod $q$ given by
 $\chi_1(g) = e^{2\pi i/(q-1)}$, to see that the set of the non trivial characters
 mod $q$ is $\{\chi_1^j \colon j=1,\dotsc,q-2\}$.
 Hence, if, for every $k\in \{0,\dotsc,q-2\}$, we denote $g^k\equiv a_k\in\{1,\dotsc,q-1\}$,
 every summation in \eqref{general-simplified}-\eqref{odd-2} and \eqref{even-simplified}
 is  of the type 
 \begin{equation}
 \label{FFT-sum}
 \sum_{k=0}^{q-2}  e\Bigl(\frac{\sigma  j k}{q-1}\Bigr) f\Bigl(\frac{a_k}{q}\Bigr),
 \end{equation}
where $e(x):=\exp(2\pi i x)$, $j\in\{1,\dotsc,q-2\}$, $\sigma=\pm 1$,
 and $f$ is a suitable function which assumes real values. 
As a consequence, such quantities are, depending on $\sigma$,
 the Discrete Fourier Transforms, or its inverse transformation, of the sequence 
 $\{ f(a_k/q)\colon k=0,\dotsc,q-2\}$. 
 This  idea was first formulated by Rader \cite{Rader1968} and it was already used in 
 \cite{FordLM2014}
 to speed-up the computation of these quantities via the use of FFT-dedicated
 software libraries. 
 
For the approach \ref{S-approach})  we can also use the \emph{decimation in frequency} 
strategy: assuming that 
in \eqref{FFT-sum} one has to distinguish between the parity of $j$
(hence on the parity of the Dirichlet characters), letting $m =(q-1)/2$,
for every $j=0,1,\dotsc,q-2$  we have that
\begin{align*}
 \sum_{k=0}^{q-2} e\Bigl(\frac{\sigma  j k}{q-1}\Bigr)  f \Bigl(\frac{a_k}{q}\Bigr)
 &
 = 
 \sum_{k=0}^{m-1}  e\Bigl(\frac{\sigma  j k}{q-1}\Bigr)  f\Bigl(\frac{a_k}{q}\Bigr)
+
 \sum_{k=0}^{m-1} e\Bigl(\frac{\sigma  j (k+m)}{q-1}\Bigr)   f \Bigl(\frac{a_{k+m}}{q}\Bigr)
  \\&
 =
 \sum_{k=0}^{m-1}  
e\Bigl(\frac{\sigma  j k}{q-1}\Bigr)
 \Bigl(
 f\Bigl(\frac{a_k}{q}\Bigr)
 +
 (-1)^{j} 
 f \Bigl(\frac{a_{k+m}}{q}\Bigr)
 \Bigr).
\end{align*}

Let now $j=2t+\ell$, where $\ell\in\{0,1\}$ and $t\in \Z$. Then, the previous equation becomes
\begin{align}
\notag
 \sum_{k=0}^{q-2}  e\Bigl(\frac{\sigma  j k}{q-1}\Bigr) f \Bigl(\frac{a_k}{q}\Bigr)
 &
 =
 \sum_{k=0}^{m-1}  
 e\Bigl(\frac{\sigma  t k}{m}\Bigr)   e\Bigl(\frac{\sigma  \ell k}{q-1}\Bigr)   
  \Bigl( 
  f\Bigl(\frac{a_k}{q}\Bigr)
 +
 (-1)^{\ell} 
 f \Bigl(\frac{a_{k+m}}{q}\Bigr)
 \Bigr)
 \\
 \label{DIF}&
 =
 \begin{cases}
 \sum\limits_{k=0}^{m-1}    e\bigl(\frac{\sigma  t k}{m}\bigr) b_k 
 & \textrm{if} \ \ell =0\\
  \sum\limits_{k=0}^{m-1}    e\bigl(\frac{\sigma  t k}{m}\bigr)  c_k  
 & \textrm{if} \ \ell =1,\\
 \end{cases}
\end{align}
 where $t=0,\dotsc, m-1$, $\sigma=\pm 1$, 
\[
b_k :=
  f\Bigl(\frac{a_k}{q}\Bigr) +  f \Bigl(\frac{a_{k+m}}{q}\Bigr)   
\quad
\textrm{and}
\quad
c_k := 
 e\Bigl(\frac{\sigma k}{q-1}\Bigr)   
 \Bigl(  f\Bigl(\frac{a_k}{q}\Bigr) -  f \Bigl(\frac{a_{k+m}}{q}\Bigr)  \Bigr).
\]
Hence, if we just need the sum  over the even, or odd, Dirichlet characters
as in the procedure  \ref{S-approach})  for $f(x)=S(x)$ or $f(x)=x$,
instead of computing an FFT transform of length $q-1$
we can evaluate  an FFT of length $(q-1)/2$, applied on a suitably
modified sequence according to \eqref{DIF}.
Clearly this represents a gain in both the speed and the memory usage
in running the actual computer program. Moreover,
if the values of $f(a_k/q)$ have to be precomputed and stored, 
this also means that the quantity of information we have to 
save during the precomputation (which will be the
most time consuming part), and to recall for
the FFT algorithm,  is reduced by a factor of $2$.

In Table \ref{table1} we give a summary of the main characteristics of both 
approaches for computing $\G_q$; it is clear that the one using $T(x)$ beats 
the one which implements $S(x)$ only in the total number of  the needed 
FFT transforms\footnote{In fact the FFT transforms can be independently
performed and hence they can be executed in parallel; this 
eliminates the unique disadvantage in using the $S$-function method.}, 
but in any other aspect  the latter is better. 
In particular the procedure \ref{S-approach})  is much faster in the 
precomputation
part since its cost is  $\le 1/2$ than approach \ref{T-approach})'s one.

\subsection{Decimation in frequency for the even Dirichlet characters case} 
 We make explicit the form that the sequence $b_k$ defined in \eqref{DIF},
 assumes in our cases.
 
It is  useful to remark   that
from $\langle g \rangle= \Z^*_q$ it trivially follows that $g^m \equiv q-1 \bmod{q}$,
where $m=(q-1)/2$.
Hence, recalling  $a_k \equiv g^k \bmod q$,   we obtain
 \(
 a_{k+m} \equiv g^{k+m}  \equiv a_k (q-1) \equiv q-a_k \bmod{q}
  \) 
  and, as a consequence, we get
\begin{equation} 
\label{f-values}
f \Bigl(\frac{a_{k+m}}{q}\Bigr)  = f\Bigl(\frac{q-a_{k}}{q}\Bigr)
=
f\Bigl(1-\frac{a_{k}}{q}\Bigr).
\end{equation}
So, inserting  the \emph{reflection formula} for $S(x)$,  see eq.~(3.3) of 
Dilcher \cite{Dilcher1994}\footnote{Pay attention to the 
fact that  the Deninger $S(x)$-function defined  in \eqref{R-def}-\eqref{S-def} is equal to 
$ -2 \log(\Gamma_1(x))$ as defined 
in Proposition 1 of Dilcher \cite{Dilcher1994}.}, into  \eqref{DIF}-\eqref{f-values},  
for every $k=0,\dotsc, m-1$  and for $f(x)=S(x)$, using \eqref{S-def}, the sequence $b_k$  becomes
\begin{align}
\notag
S\Bigl(\frac{a_k}{q}\Bigr) &+  S \Bigl(\frac{a_{k+m}}{q}\Bigr)   
=
S\Bigl(\frac{a_k}{q}\Bigr) +  S\Bigl(1-\frac{a_{k}}{q}\Bigr)
=   
\\
\label{S-DIF-formula}
&=  
\Bigl(\log \frac{a_k}{q} \Bigr)^2 + 
\sum_{n=1}^{+\infty} 
\Bigl(
\Bigl(\log \bigl(n+\frac{a_k}{q} \bigr)  \Bigr)^2  
+
\Bigl(\log \bigl(n-\frac{a_k}{q} \bigr) \Bigr)^2
- 
2 (\log n)^2
\Bigr),
\end{align}
where   $a_k \equiv g^k \bmod q$, while, using \eqref{S-alt-def}, we obtain
\begin{align}
\notag
S\Bigl(\frac{a_k}{q}\Bigr) +  S\Bigl(1-\frac{a_{k}}{q}\Bigr)
&
=    
2 \int_{0}^{+\infty} 
\Bigl(- e^{-t} + \frac{e^{-\frac{a_{k}}{q}t} +e^{-(1-\frac{a_{k}}{q})t} - 2e^{-t}}{1-e^{-t}} \Bigr)
\frac{\gamma + \log t}{t} \ \dx t
\\&
 \label{S-alt-DIF-formula}
 =
2 \int_{0}^{+\infty} 
\Bigl(- 3 + e^{-t} +e^{\frac{a_{k}}{q}t} +e^{(1-\frac{a_{k}}{q})t}  \Bigr)
\frac{\gamma + \log t}{t(e^t-1)} \ \dx t ,
\end{align}
in which we exploited the uniform convergence of the involved integrals.
To optimise speed and precision,
both equations \eqref{S-DIF-formula}-\eqref{S-alt-DIF-formula} will be used during the 
actual computations; when possible we will exploit the exponential decay $e^{-ct}$,
with   $c=\min (a_k/q,1-a_k/q)$, of the
integrand function in \eqref{S-alt-DIF-formula}  using the PARI/GP function
 {\tt intnum}. But when the parameter $c$ will become too small to give reliable results, 
 we will switch to apply the PARI/GP function {\tt sumnum} to eq.~\eqref{S-DIF-formula};
  in this case, roughly speaking,  the decay
 order is  $ (\log n) /n^2$.
 
Hence, thanks to the previous formulae,  the 
number of calls to the {\tt sumnum} or  {\tt intnum}  functions
 required  in the precomputation of the $S$-values is
 reduced by a factor of $2$ with respect to the ones needed 
 to  precompute  the $T$-values.

If  we are just interested in   the computation of $\G_q^+$,   we can directly
 use \eqref{EK+-formula} 
in which we can embed \eqref{S-DIF-formula}-\eqref{S-alt-DIF-formula} and the following remark about 
 the needed $\log \Gamma$-values. Assuming $f(x) = \log \Gamma(x)$, 
 using \eqref{f-values}
 and  the well-known \emph{reflection formula} $\Gamma(x) \Gamma(1-x)  = \pi / \sin(\pi x)$,
we obtain
\begin{align*}
\log \Bigl(\Gamma \bigl(\frac{a_k}{q}\bigr)\Bigr) 
+ \log \Bigl(\Gamma \bigl(\frac{a_{k+m}}{q}\bigr)\Bigr)
=
\log \Bigl(\Gamma \bigl(\frac{a_k}{q}\bigr)\Bigr) 
+ \log \Bigl(\Gamma \bigl(1-\frac{a_{k}}{q}\bigr)\Bigr)
=  
 \log \pi - \log\Bigl(\sin\bigl( \frac{\pi a_k}{q}\bigr)\Bigr),
\end{align*}
thus further simplifying the final computation  by replacing the $\Gamma$-function 
with the $\sin$-function.

\subsection{Decimation in frequency for the odd Dirichlet characters case}
\label{DIF-odd}
  We make explicit the form that the sequence $c_k$ defined in \eqref{DIF},
 assumes in our cases.
 
If  we are just interested in   the computation of $\G_q-\G_q^+$,  we can directly
use \eqref{diff-equation};  using the   \emph{reflection formula} 
$\Gamma(x) \Gamma(1-x)  = \pi / \sin(\pi x)$
and arguing as in  the previous paragraph,
we obtain
\begin{align*}
\log \Bigl(\Gamma \bigl(\frac{a_k}{q}\bigr)\Bigr)  
- \log \Bigl(\Gamma \bigl(1-\frac{a_{k}}{q}\bigr)\Bigr)
&=
2 \log \Bigl(\Gamma \bigl(\frac{a_k}{q}\bigr)\Bigr)
+
\log\Bigl(\sin\bigl( \frac{\pi a_k}{q}\bigr)\Bigr) - \log \pi,
\end{align*}
for every $k=0,\dotsc, m-1$, $m=(q-1)/2$,
and hence $c_k$ is modified accordingly.  In this case the gain of using the previous formula is
that the number of needed evaluations of the $\log\Gamma$-function is reduced by a factor of $2$.
 
 The case  in which $f(x)=x$ is easier; using again $\langle g \rangle = \Z^*_q$, 
 $a_k \equiv g^k \bmod q$ and   $g^m \equiv q-1 \bmod{q}$,
 we can write that $ a_{k+m}  \equiv  q-a_{k} \bmod{q}$;  hence
\[
 a_k  -   a_{k+m} =
 a_k -(q-a_{k})
=
 2a_k- q,
\]
so that in this case we obtain $c_k=   e(\sigma k / (q-1))(2a_k/q -1)$ for every 
$k=0,\dotsc m-1$, $m=(q-1)/2$, $\sigma=\pm 1$.

\subsection{Computations trivially summing over $a$ (slower but with more digits available).}
Unfortunately in \texttt{libpari} the FFT-functions  work only if $q=2^\ell+1$, for some $\ell\in \N$. 
So  we had to trivially perform these summations and  hence, in practice,  
this part is the most time consuming one in both the procedures  
\ref{T-approach}) and \ref{S-approach}) since it has a quadratic cost 
in $q$.
Being aware of such limitations, we used PARI/GP (with the trivial way to compute 
the sum over $a$) 
to evaluate $\G_q$ and $\G_q^+$ with  these three approaches for every  odd prime $q\le 300$, 
on a   Dell OptiPlex-3050 (Intel i5-7500 processor, 3.40GHz, 
16 GB of RAM and running Ubuntu 18.04.2) using a   precision of $30$ digits, see Table \ref{table2};
we also inserted there the values of $M_q$, defined in \eqref{Mq-def}, for the same set of primes. 
Such results largely extend the precision of the data in Table 1 on page 1472  of \cite{FordLM2014}. 
The computation of the values of Table \ref{table2} needed  19 seconds
 using the $S$-function, 33 seconds   using the $T$-function
 and  51 seconds using PARI/GP \texttt{lfun} function.
We also computed the values of $\G_q$ and $\G_q^+$,  with a precision of 30  digits, for 
$q=1009,2003,3001, 4001, 5003,$ $6007, 7001,$ $8009,$ $9001,10007$, $20011, 30011$,
as you can see in Table \ref{table3}.
These numbers were chosen to heuristically evaluate how the computational cost depends
on the size of $q$.
In this case, in the fifth column of Table \ref{table3} we also reported the running time of 
the direct approach,
\emph{i.e.} using \eqref{EKq-def}, the third and fourth columns are respectively the 
running times
of the other two procedures.
For these values of $q$ it became clear that the computation time spent in performing 
the sums 
over $a$  was the longest one. This means that inserting an FFT-algorithm  is fundamental to further 
 improve the performances of both the algorithms \ref{T-approach})-\ref{S-approach}).
We discuss this in more detail in the next paragraph.

\subsection{Computations  summing over $a$ via FFT (much faster but with less 
digits available).} 
As we saw before,  for large $q$  the time spent in summing over $a$
dominates the computational cost.  So we implemented the use of FFT for this task.
We  first used the {\tt gp2c} compiler tool to 
obtain suitable C programs to perform the precomputations of the 
needed $T$ and $S$-values with $38$ digits
and save them to the hard disk\footnote{If we do not use the GSL to directly compute
$\psi$, we need
to insert its precomputation here.}. Then  we passed 
such values to the C programs which used the 
{\tt fftw} \cite{FFTW} software library  to perform the  FFT step. 
In such a final stage the performance was thousand times faster
than  the one for the same stage trivially performed; 
as an example you can compare the running
times for $q=10007$, $20011$, $30011$ in Tables \ref{table3} and \ref{table4}. 
The running times for the approaches   \ref{T-approach}) and \ref{S-approach})
reveal that the latter is faster, mainly because it requires less input operations to 
gain the stored
precomputed information since the FFT works on a set of data of half the length 
than in the former case\footnote{If $\psi$
is precomputed using PARI/GP, then the gain ratio in the stored space and in the
number of input/output operations is raised to $3/4$.}.

This way we computed the values of $\G_q$ and $\G_q^+$
for $q=40009$, $42611$, $50021$, $60013$, $70001$,  $80021$,  $90001$, $100003$, $305741$, $1000003$, 
$4178771$, $6766811$, $10000019$, $28227761$, $75743411$
with the long double precision, see Table \ref{table4}. 
These computations were performed 
with the Dell OptiPlex machine mentioned before.

 Some of these $q$-values were chosen for their size and others with the help of
  $\B$, the ``greedy sequence of prime offsets'', 
 \url{http://oeis.org/A135311}, in the following way.
 We define $\B$ using induction, by 
   $b(1)=0\in \B$ and $b(n)\in \B$ if it is  the smallest integer exceeding $b(n-1)$ such that 
 for every prime $r$ the set $\{b(i) \bmod r\colon 1\le i \le n\}$ has at most $r-1$ elements. 
 An equivalent statement, assuming that the prime $k$-tuples conjecture holds,
 is that $b(n)$ is minimal such that $b(1)=0$ and  there are infinitely many primes $q$ with $b(i)q+1$ prime
 for $2\le i \le n$, $n\ge 2$.   
Let now
\[m(\A):= \sum_{i=1}^s\frac{1}{a_i},\]
 where $\A$ is an admissible set, \emph{i.e.}, $\A=\{a_1,\dotsc, a_s \}$, $a_i\in \N$, $a_i\ge 1$, such that
 does not exist a prime $p$ such that $p\mid n\prod_{i=1}^s (a_in+1)$ for every $n\ge 1$.
Thanks to Theorem 2 of
Moree \cite{Moree2018}, if the prime $k$-tuples conjecture holds and if $\A$ is an admissible set, then
$\G_q < (2- m(\A)+ \odi{1})  \log q$ for $\gg x/(\log x)^{-\vert \A \vert -1}$ primes $q \le x$.
Moreover, by Theorem 6 of Moree \cite{Moree2018}, assuming both the Elliott-Halberstam 
and the prime $k$-tuples conjectures,  if $\A$ is an admissible set then
$\G_q = (1- m(\A)+ \odi{1})  \log q$ for $\gg x/(\log x)^{-\vert \A \vert -1}$ primes $q \le x$.

 The greedy sequence of prime offsets $\B$ has the property that any finite subsequence is an admissible set.
 With a PARI/GP script we   computed the first $2089$ 
elements of $\B$  since for $\CC := \{b(2), \dotsc, b(2089)\}$ we get  $m(\CC) >2$.
So, if we are looking for negative values of $\G_q$, it seems to be a good criterion to evaluate $\G_q$ for 
a  prime number $q$  such that $bq+1$ is prime for many elements $b\in \CC$ 
(clearly it is better to start with the smaller available $b$'s).
 To be able to measure this  fact, we define
\begin{equation}
\label{vq-def}
v(q) :=  \sum_{\substack{2\le i\le 2089; \ b(i)\in \CC\\ b(i)q+1\ \text{is prime} }}\frac{1}{b(i)}.
\end{equation}
Some of the $q$-values  written before in this paragraph are such that $v(q)>1.15$ so that, 
thanks to Moree's results
already cited,  they are good candidates  to have a negative Euler-Kronecker constant.
 The complete list of $q\le 10^{10}$ such that $v(q) >1.2$
 is towards the end of the PARI/GP script {\tt testseq}
 that can be downloaded here:
\url{http://www.math.unipd.it/~languasc/EK-comput.html}.

\subsubsection{Data for the scatter plots.}
\label{data-plots}
After having evaluated the running times of the previous examples, 
we decided to provide the colored
scatter plots, see Figures \ref{fig1}-\ref{fig2}, of the normalised values 
of  $\G_q$ and $\G_q^+$ (both  in long double precision) for every odd prime $q\le 10^6$ 
thus  enlarging the known range of the data on $\G_q$
and $\G_q^+$, see \cite{FordLM2014}.
For performing the needed precomputations of the $S$-values, we used
the cluster of the  Department  of Mathematics
of the University of Padova; the cluster  setting is described here: 
\url{http://computing.math.unipd.it/highpc}. 
The  minimal value of $\G_q/\log q$, $3\le q\le 10^6$, $q$ prime,  is $0.13067\dotsc$ and 
it is attained at $q=305741$, as expected;   the  maximal value    is $1.62693\dotsc$ and 
it is attained at $q=19$.
 The  minimal value of $\G^+_q/\log q$, $3\le q\le 10^6$, $q$ prime,  is $0.451468\dotsc$ 
 and it is attained at $q=918787$;  the  maximal value    is $1.42626\dotsc$ and it is attained 
 at $q=2053$. 
 The points $(q,\G_q/\log q)$ and $(q,\G^+_q/\log q)$ in Figures \ref{fig1}-\ref{fig2}  
 are  
 colored 
 in orange  if $v(q) \le 0.25$ ($65.65$\% of the cases), 
 in green if $0.25<v(q) \le 0.5$ ($23.62$\%), 
 in blue if $0.5<v(q) \le 0.75$ ($6.29$\%), 
 in black if $0.75<v(q) \le 1$ ($4.21$\%), 
 and in red if $v(q) >1$ ($0.23$\%). 
 The behaviour of $\G_q$ is the expected one since the red strip 
 essentially corresponds with its minimal values, while the minima of $\G^+_q$  seem to be 
 less related to $v(q)$; we plan to investigate this phenomenon in the next future.
 The complete list of numerical results for $\G_q$  and $\G^+_q$ can be downloaded
 at the following web address:
\url{https://www.math.unipd.it/~languasc/EKcomput/results}.

\subsubsection{Computations for larger $q$.}
For  values  of $q$ larger than $30$ millions  the precomputation of $T$ and $S$, if performed 
on a single desktop computer,  would require too much  time; hence
we parallelised them on the cluster  previously mentioned. To   check the correctness of such  
computations
it is possible to use the following formulae; recalling that
$\gamma = 0.577215664901\dotsc$ and 
$\zeta^{\prime\prime}(0) =  - 2.006356455908\dotsc$, 
we have that 
 \begin{align}
 \label{S-sum}
 \sum_{a=1}^{q-1} S\Bigl(\frac{a}{q}\Bigr) & = -\zeta^\second(0)(q-1) - \log q \log(2\pi) - \frac{(\log q)^2}{2}  ,
\\
 \label{psi-T-sum} 
 \sum_{a=1}^{q-1} T\Bigl(\frac{a}{q}\Bigr) &= \frac{q}2  (\log q)^2 + \gamma q \log q .
 \end{align}
Formula \eqref{S-sum} is an immediate consequence
of Theorem 2.5 of Deninger \cite{Deninger1984}  and formula  
\eqref{psi-T-sum} follows   from equation  (7.10) of Dilcher \cite{Dilcher1992}.

Moreover, for being able to handle very large cases, we  used a dedicated 
{\tt fftw} interface\footnote{It is called the {\tt guru64}
 interface; see the user's manual of {\tt fftw} \cite{FFTW}.}
 which is able to perform transforms whose length is greater than $2^{31}-1$.
 
In this way we were able to  obtain an  independent  confirmation 
of Theorem 4 of \cite{FordLM2014}  
getting 
$\G_{964477901}=  -0.18237472563711916085\dotsc$, 
since we computed it using the quadruple precision.
At the same time we also got $\G^+_{964477901} = 10.40222338242826353694 \dotsc$ 
To do so we first split the computation, 
 with a precision of 38 digits,
 of the needed decimated in  frequency values of  $S$
in $49$ subintervals $I_j$ of size $10^7$ each (for $T$ we would need $97$ intervals of such a length);  
the computation time required for each $I_j$ was on average about $1600$  
minutes  on one of the cluster's machines.
Then we passed such values to the programs that performed 
the FFT-step and got the final results. This last part needed  about 23 minutes  (long double precision)
or  522 minutes  (quadruple precision) of computation time
on an Intel(R) Xeon(R) CPU E5-2650 v3 @ 2.30GHz, with 160 GB of RAM, running Ubuntu 16.04.
A similar procedure let us to get analogous  computation times  for  the long double precision evaluation of
$\G_{1217434451}=  0.877596 \dotsc$ and
$\G^+_{1217434451} =  12.946690 \dotsc$  
 
We then looked for prime numbers $q$ such that $v(q)>v(964477901)=1.2369344\dotsc$
and we found that   $v(2918643191) = 1.2440460\dotsc$
In about 90  minutes of computation time for the FFT-step on the same 
 machine mentioned before we  got that 
$\G_{2918643191}=  0.302789\dotsc$
and $\G^+_{2918643191} = 12.573983\dotsc$,
using the long double precision.
In this case it seems that procedure in \ref{T-approach})
is much less stable than the one in \ref{S-approach})
probably because of the fact that 
$T(x)$ and $\psi(x)$ are much larger, for $x\to 0^+$, than $S(x)$ and $\log(\Gamma(x))$.
Computations for further  ``good'' candidates, in the sense that $v(q)>1.18$,  like 
$q =193894451$, $212634221$, $251160191$
$538906601$, $1139803271$, $1217434451$, $1806830951$, 
$2488788101$,
$2830676081$, $7079770931$
were also performed. 
The computations for these primes were performed on the cluster previously mentioned.

Moreover, for $q=9109334831$ we got that
$\G_{9109334831}=  -0.248739\dotsc$,  thus obtaining a new minimal value for $\G_q$
 and a new example   
of Theorem 4 of \cite{FordLM2014}; at the same time we also got
$\G^+_{9109334831} = 12.128187\dotsc$ The precomputations for this case,
performed with the same strategy used for the smaller primes $q$ mentioned in this paragraph,
required about nine days on the cluster and the FFTs computation required about  1000 minutes
on the Xeon machine mentioned before (this amount of time also depends on a runtime RAM swapping phenomenon) 
or 312 minutes on the new CAPRI infrastructure of the University of Padova
(``Calcolo ad Alte Prestazioni per la Ricerca e l'Innovazione''; whose CPU is  an
Intel(R) Xeon(R) Gold 6130 CPU @ 2.10GHz, with 256 cores and equipped with 6TB of RAM).
Such a result was then double-checked on CAPRI  using the much slower algorithm \ref{T-approach}).
A further new example of Theorem 4 of \cite{FordLM2014} we obtained is 
$\G_{9854964401}=  -0.096465\dotsc$ which required about ten days of time for the precomputations
and 326 minutes for the FFT stage on CAPRI.
As usual the result was double-checked  using the approach \ref{T-approach}).

All the results  mentioned in this paragraph
are collected in Table \ref{table5}.
The PARI/GP scripts and the C programs used and the computational results obtained
are available at the following web address:
\url{http://www.math.unipd.it/~languasc/EK-comput.html}.

\section{On the absolute value of the logarithmic derivative\\of Dirichlet $L$-functions}
\label{Mq-eval}
Using \eqref{derlog-01-rel}-\eqref{L0}, \eqref{even-1} and \eqref{R-S-summed}-\eqref{even-2}, for every
odd prime $q$ we immediately get 
\begin{equation*} 
M^{\textrm{odd}}_q
:=\max_{\chi\, \textrm{odd}} \
\Bigl\vert\frac{L^\prime}{L}(1,\chi) \Bigr\vert  
= 
\max_{\chi\, \textrm{odd}} \
 \Bigl\vert
\gamma + \log(2\pi)  
 +
\frac{1}{B_{1,\overline{\chi}} }
 \sum_{a=1}^{q-1}   \overline{\chi}(a)\log\Bigl(\Gamma\bigl(\frac{a}{q}\bigr)\Bigr)
 \Bigr\vert 
\end{equation*}  
and
\begin{equation*} 
M^{\textrm{even}}_q
:=\max_{\substack{\chi \neq \chi_0\\ \chi\, \textrm{even}}} \
\Bigl\vert\frac{L^\prime}{L}(1,\chi) \Bigr\vert  
=
\max_{\substack{\chi \neq \chi_0\\ \chi\, \textrm{even}}} \
 \Bigl\vert
\gamma + \log(2\pi) 
-\frac{1}{2} 
\frac{\sum_{a=1}^{q-1} \overline{\chi}(a) \ S(a/q)}
{\sum_{a=1}^{q-1} \overline{\chi}(a)\log\bigl(\Gamma(a/q)\bigr)}
 \Bigr\vert.
\end{equation*}
Hence we can compute  $M_q=\max_{\chi \neq \chi_0} \vert L^\prime/L(1,\chi) \vert  
= \max(M^{\textrm{odd}}_q, M^{\textrm{even}}_q)$
 using the  values of $\log\Gamma$ and $S$ 
obtained for the computation of $\G_{q}$ and $\G^+_{q}$.
In Table \ref{table2} we give the values of $M_q$  for every odd prime up to $300$  computed, 
using PARI/GP, with 
a precision of $30$ digits.
Using the data  in Section \ref{data-plots} we also computed, on the Dell Optiplex machine 
previously mentioned, 
the values of $M_q$ and $M_q/\log\log q$ for every odd prime $q\le 10^6$ 
and in Figures \ref{fig3}-\ref{fig4}  we inserted  their scatter plots that largely extend Figure 1 
of  Ihara-Murty-Shimura \cite{IharaMS2009} (please remark that our $M_q$ is denoted as  
$Q_m$ there).
Such data in Figures \ref{fig3}-\ref{fig4} also
fit, for $q$ sufficiently large, with  the estimate $M_q \leq (2+\odi{1}) \log \log q$ as 
$q$ tends to infinity, proved,
under the assumption of GRH, in
Theorem 3 of \cite{IharaMS2009}.
We also remark  that $M^{\textrm{odd}}_q > M^{\textrm{even}}_q$ for $62521$ 
cases  over  a total number of primes
equal to $78497$ ($79.65\%$)
and that $M^{\textrm{even}}_q > M^{\textrm{odd}}_q$ in the remaining $15976$ cases  ($20.35\%$).

 The complete list of numerical results for $M_q$  
 can be downloaded
 at the following web address:
\url{https://www.math.unipd.it/~languasc/EKcomput/results}. 

\section{On the generalised Euler constants in arithmetic progressions $\gamma_{k}(a,q)$} 
\label{gen-euler-arit-prog}
Recall that $q$ is an odd prime.
In the case   we have to precompute 
$T(a/q)$ and we also need $\psi(a/q)$.
Hence, as a by-product we can also obtain the values of the generalised Euler
constants $\gamma_{0}(a,q)$ and $\gamma_{1}(a,q)$, 
see subsections \ref{gamma0-aq-comp}-\ref{gamma1-aq-comp}. In practice this 
is done by activating an optional flag in the main gp script.
The computation of  $\gamma_{k}(a,q)$ for $k\ge 2$ is described in subsection \ref{gammak-aq-comp}.
\subsection{Generalised Euler constants $\gamma_{0}(a,q)$} 
\label{gamma0-aq-comp}
For 
$\gamma_{0}(a,q)$ with  $1\leq a \leq q-1$, $q$  odd prime, by \eqref{gammakaq-def}
we have
\[
\gamma_{0}(a,q) 
=
- \frac{1}{q}
\Bigl(
\log q
+  \psi(\frac{a}{q})
\Bigr). 
\]
Recalling that $\psi(1)=-\gamma$, we also have
\(
\gamma_{0}(q,q)
=(\gamma -  \log q)/q.
\)

\subsection{Generalised Euler constants $\gamma_{1}(a,q)$}
\label{gamma1-aq-comp}
For 
$\gamma_{1}(a,q)$ with $1\leq a \leq q-1$,  $q$  odd prime, we can use  \eqref{gammakaq-def}
and \eqref{T-def}. 
This way we get
\begin{align*}
\gamma_{1}(a,q) 
&=
- \frac{1}{q}
\Bigl(
\frac{(\log q)^{2}}{2}
+ (\log q) \psi(\frac{a}{q})
+
\psi_1(\frac{a}{q}) 
\Bigr)
=
\frac{1}{q}\
\Bigl(\gamma_1-\frac{(\log q)^{2}}{2}
-(\log q) \psi(\frac{a}{q})
-
T\bigl(\frac{a}{q}\bigr)
\Bigr).
\end{align*}
Moreover, since $\psi(1)= -\gamma$ and $T(1)=0$, we also have
\[
\gamma_{1}(q,q)
=\frac{1}{q}\Bigl(\gamma_1+\gamma \log q -\frac{(\log q)^{2}}{2}\Bigr).
\]

Using the formulae in the previous two paragraphs we computed
 $\gamma_{0}(a,q)$ and $\gamma_{1}(a,q)$ with $q$ prime, $3\le q \le 100$,
 $1\le a\le q$, in about 4  seconds of computation time with a precision of 30 digits.

 Such results are listed towards the end of the gp-script file that can be downloaded here: 
 \url{http://www.math.unipd.it/~languasc/EK-comput.html}.

 \subsection{The general case $\gamma_{k}(a,q)$,  $k\ge 2$}
 \label{gammak-aq-comp}
% \label{Stieltjes} 
 
The general  case $\gamma_{k}(a,q)$, $k\in \N$,  $k\ge 2$, $q\ge 1$,
 $1\leq a \leq q$,     do not follow from the data already computed for the Euler-Kronecker constants
 since we need information about the values of $\psi_n(x)$, for every $2\le n \le k$.   
Such a direct computation of both $ \psi_n(a/q)$
and $\gamma_n$ can be easily performed via eq.~\eqref{gammakaq-def}-\eqref{gen-digamma-def} 
using the PARI/GP summing function
 {\tt sumnum} paying attention to submit a sufficiently fast convergent 
sum. For example, to compute $\gamma_n$, $n\in \N$, we  used the formulae
\begin{equation}
\label{gen-euler-const-comp}
\gamma_n   =
 \sum_{m=1}^{+\infty}
   \Bigl( 
   \frac{(\log m)^n}{m} -   \frac{1}{n+1}
   \sum_{j=0}^n\binom{n+1}{j} (\log m)^{j} (\log \bigl(1+\frac{1}{m})\bigr)^{n+1-j} 
   \Bigr)
\end{equation}
and
\begin{equation}
\label{gen-euler-const-comp1}
\gamma_n   =
 \sum_{m=1}^{+\infty}
  \Bigl( 
(\log m)^n   \Bigl( 
   \frac{1}{m} - \log\bigl(1+\frac{1}{m}\bigr)\Bigr)
   -
    \frac{1}{n+1}
   \sum_{j=0}^{n-1}\binom{n+1}{j} (\log m)^{j} (\log \bigl(1+\frac{1}{m})\bigr)^{n+1-j} 
   \Bigr),
\end{equation}
which both easily follow  from \eqref{gamman-def}.
We get,  in less than  7 seconds of  time and 
with a precision of at least $40$ digits, the  results in Table \ref{table6}; 
to be sure about the correctness of such results
we computed them twice  using the formulae 
\eqref{gen-euler-const-comp}-\eqref{gen-euler-const-comp1}
and then we compared the outcomes.
These values are in agreement with the data on page 282 of Bohman-Fr\"oberg \cite{BohmanF1988}
 for $n=0,\dotsc,20$.
For larger $n$'s the formulae in 
\eqref{gen-euler-const-comp}-\eqref{gen-euler-const-comp1}
seem to be not good enough
to get precise results via the {\tt sumnum} function with this precision level.

To compute $\psi_n(a/q)$ and, as a consequence,  $\gamma_k(a,q)$,  we can proceed in a similar
way as we did for $T(a/q)$
and $ \gamma_1(a,q)$,
see the program  {\tt Gen-Euler-constants.gp}  here \url{http://www.math.unipd.it/~languasc/EK-comput.html}.
 Towards the end of this  program   file you can find a large list
(too long to be included here) of computed values of $\gamma_k(a,q)$
for $1\le k \le 20$, $1\le q \le 9$, $1\le a \le q$,  with a precision of $20$ digits.
In about 50 seconds of computation time we  replicated Dilcher's computations, 
since the values we got  are in agreement with the data on pages S21-S24 of \cite{Dilcher1992}.

\renewcommand{\bibliofont}{\normalsize}

 \vskip1cm
\noindent 
Alessandro Languasco,     
Universit\`a di Padova,
 Dipartimento di Matematica,
 ``Tullio Levi-Civita'',
Via Trieste 63,
35121 Padova, Italy.    
{\it e-mail}: alessandro.languasco@unipd.it 

 \newpage
\section{Tables}
\label{tables}
\newcommand{\footnoteref}[1]{\textsuperscript{\ref{#1}}}
\renewcommand{\thefootnote}{\fnsymbol{footnote}} 

\begin{table}[htp]
\scalebox{0.75}{
\rotatebox{90}{
\renewcommand{\arraystretch}{1.8}%
\begin{tabular}[h]{|l|c|c|c|}   
\hline
 Comparison &   \textbf{Approach \ref{T-approach})} &  \textbf{Approach \ref{T-approach})} & \textbf{Approach \ref{S-approach})} \\  
  &   ($\psi$ comp. with GSL)   & ($\psi$ comp. with PARI/GP)  &    \\ \hline
   \hline
\textbf{Magnitude of the functions for $x\to 0^+$}:    &     $\psi(x)\sim  -1/x$&   $\psi(x)\sim  -1/x $&$\log\bigl(\Gamma(x)\bigr)\sim \log(1/x)$  \\   
 &   $T(x)\sim \frac{\log(1/x)}{x}$   & $T(x)\sim \frac{\log(1/x)}{x}$&  $S(x) \sim   (\log x)^2$   \\ \hline    
\hline  
\textbf{Precomputations ($T$ and $S$ with PARI/GP)}:  &     &   &\\ \hline  
needed space for storing precomputed values   &   $q-1$   values of  &  $2(q-1)$   values of  & $(q-1)/2$ values of    \\   
 ($\langle g\rangle = \Z_q^*$, $a_k:=g^k \bmod q$):   &     $T(a_k/q)$  & $T(a_k/q)$ and $\psi(a_k/q)$& $S(a_k/q)+S(1-a_{k}/q)$  \\ \hline  number of {\tt write} operations on  hard disks:&   $q-1$   &$2(q-1)$ & $(q-1)/2$   \\ \hline  
number of {\tt sumnum} or  {\tt intnum} calls: &   $q-1$   & $q-1$ & $(q-1)/2$   \\ \hline  
\hline 
\textbf{FFT-step (with {\tt fftw})}:  &     &  &\\ \hline 
number of  {\tt read} operations on  hard disks:&   $q-1$  & $2(q-1)$ &  $(q-1)/2$   \\ \hline  
number of FFTs:  &  $2$  & 2 & $3$ \\ \hline  
length of FFTs:   &  both $q-1$  &  both $q-1$ & one of length $q-1$;   \\   
    &    &  & the others of length $ (q-1)/2$      \\ \hline  
total RAM usage (in number of long &    & &  \\   
 double positions;  in-place FFTs): &  $2q+2$  & $2q+2$& $2q$\footnote[2]{But  the computation for $\G_q^+$ requires  only $(3q+5)/2$ long double positions; so, reusing  a portion of the RAM after the computation of $\G_q-\G_q^+$,
 in the second part of the program we essentially have a gain of about $(q-1)/ 2$  long double positions for the RAM usage.} \\
\hline  
\end{tabular}  
}
}
\caption{\label{table1}
Comparison of  the main characteristics of procedures \ref{T-approach}) and \ref{S-approach})
to compute $\G_q$ and $\G_q^+$.}
\end{table} 

\begin{table}[htp]
\scalebox{0.7}{
\begin{tabular}[h]{|c|c|c|c|}
\hline
$q$  &  $\G_{q}$ &  $\G_{q}^+$ & $M_q$\\ \hline
$3$ & $0.94549728087168070323974999415\dotsc$ & $ 0.57721566490153286060651209008\dotsc$ & $0.36828161597014784263323790407\dotsc$  \\
$5$ & $1.72062421251340476169572878865\dotsc$ & $1.40489514161703774859755907976\dotsc$ & $0.82767947671550488799104698967\dotsc$  \\
$7$ & $2.08759407471733013281542471957\dotsc$ & $1.95715645444971475271382186143\dotsc$ & $0.69374325299917902224231637393\dotsc$  \\
$11$ & $2.41542590428326783034287963583\dotsc$ & $2.66207409890433174906654072453\dotsc$ & $0.64960999942397995363690453077\dotsc$  \\
$13$ & $2.61075773741765019699776108857\dotsc$ & $2.89959572414790509559591203013\dotsc$ & $0.69630986299203715584089218352\dotsc$  \\
$17$ & $3.58197604409757765927178812919\dotsc$ & $3.23179164885108167689200470642\dotsc$ & $1.36293176857311326439833395890\dotsc$  \\
$19$ & $4.79040941571428332590703936458\dotsc$ & $3.36702810226943360422911738361\dotsc$ & $1.56821936415476775304938942269 \dotsc$  \\
$23$ & $2.61128917618820092550739164964\dotsc$ & $3.56605274186303485506490005633\dotsc$ & $1.07370241439895666993863022504\dotsc$  \\
$29$ & $3.09373170599426872316275179819\dotsc$ & $3.77451272291818155837540505527\dotsc$ & $1.37173438584080190328583030799\dotsc$  \\
$31$ & $4.31444292526747509770757441042\dotsc$ & $3.74063417131631765163927862231\dotsc$ & $ 1.41315141911004437078399808370\dotsc$  \\
$37$ & $4.30493818995760201798557926417\dotsc$ & $3.88346103237113739135523493388\dotsc$ & $1.29518958101078356915278401821\dotsc$  \\
$41$ & $3.97152162792133216028257040014\dotsc$ & $3.90067243331576039538420460289\dotsc$ & $1.29673609198958173353796568380\dotsc$  \\
$43$ & $4.37862750574695049413775062336\dotsc$ & $4.37462848511375110150884874389\dotsc$ & $1.41176882240051173489451389181\dotsc$  \\
$47$ & $4.79939425890741613452758429988\dotsc$ & $4.78330592374031492736088514964\dotsc$ & $1.39567565425273602292102717603\dotsc$  \\
$53$ & $4.33773685859709231869696082307\dotsc$ & $4.06734814093911422415451881781\dotsc$ & $1.30627572903790815149667975264\dotsc$  \\
$59$ & $5.43351634538500398077634438193\dotsc$ & $5.74977495098717868985714511291\dotsc$ & $1.81899383678937843989348366929\dotsc$  \\
$61$ & $5.07108519057651619595805098113\dotsc$ & $4.71919160448137601223479232791\dotsc$ & $1.41809980889441627035459190983\dotsc$  \\
$67$ & $5.29213930662896260873428461831\dotsc$ & $5.49478574409231087894450914285\dotsc$ & $1.67019193303154369921782607634\dotsc$  \\
$71$ & $5.25525819281894616772013128637\dotsc$ & $5.02459221437013823603453457463\dotsc$ & $1.47455511100236771011015896767\dotsc$  \\
$73$ & $4.06694909044749529201648815625\dotsc$ & $5.56638018904420607773144876527\dotsc$ & $1.78248970799598673447282517891\dotsc$  \\
$79$ & $4.99827631817068010789431392945\dotsc$ & $4.31392816983842153234814442952\dotsc$ & $1.34616837027813468918588610688\dotsc$  \\
$83$ & $3.03313611343607418716403819105\dotsc$ & $4.06119890648015486954960478374\dotsc$ & $1.34527786237910789501875868023\dotsc$  \\
$89$ & $4.16409079888983276880841110372\dotsc$ & $5.44834851555434719261902953243\dotsc$ & $1.61654649274126300156782088673\dotsc$  \\
$97$ & $4.89124074040389666830751468857\dotsc$ & $4.44563411256346738186380452664\dotsc$ & $1.60286118570076458480362218799\dotsc$  \\
$101$ & $5.29701289150966971887860032739\dotsc$ & $5.93364557387726998305789899164\dotsc$ & $1.51871979857079618912367283335\dotsc$  \\
$103$ & $5.14433955125208822113330503220\dotsc$ & $5.53312508630999898815400644939\dotsc$ & $1.56072764165486011343921965820\dotsc$  \\
$107$ & $5.45827420997024503421680245453\dotsc$ & $5.35744691959596839332603590620\dotsc$ & $1.55529418086936504978552066530\dotsc$  \\
$109$ & $6.90663814626423653219469837704\dotsc$ & $6.28639312060842026587282318484\dotsc$ & $1.65357828827908326582841136643\dotsc$  \\
$113$ & $4.02173038257803067578318006617\dotsc$ & $4.71308052553071355344451609738\dotsc$ & $1.51486982889352164427060492878\dotsc$  \\
$127$ & $5.08859912415333449423215636240\dotsc$ & $5.28427526641642291108714895825\dotsc$ & $1.55590143040596443193792941854\dotsc$  \\
$131$ & $2.83682634158837909860285797321\dotsc$ & $4.29182422162389365669036230041\dotsc$ & $1.43797882292531602089564238879\dotsc$  \\
$137$ & $4.93700022614368468691962999711\dotsc$ & $5.17281966401368126952267004684\dotsc$ & $1.53929870904867707257469538680\dotsc$  \\
$139$ & $5.88916863399867186726383730369\dotsc$ & $5.15673467267785693456200640445\dotsc$ & $1.58828875478913218915240825692\dotsc$  \\
$149$ & $5.98342477769515981450242785739\dotsc$ & $6.35744273145487616682151978517\dotsc$ & $1.55933423387754689170927007457\dotsc$  \\
$151$ & $5.04201611352872179914519461022\dotsc$ & $5.66732269410388218441768644382\dotsc$ & $1.48171078244888795642226012230\dotsc$  \\
$157$ & $7.40802206572222729350845201390\dotsc$ & $5.67766459100970078752076942990\dotsc$ & $ 1.52915091159611605159149879696\dotsc$  \\
$163$ & $5.92966482288720678755499913844\dotsc$ & $5.54289611872522541669860167904\dotsc$ & $2.16832712928352380386400324642\dotsc$  \\
$167$ & $8.03300175268872470467583357802\dotsc$ & $6.80394798958259907108839110755\dotsc$ & $1.56607236656750344030293511154\dotsc$  \\
$173$ & $3.38434753653206190344297798897\dotsc$ & $4.74313680866654143318864467269\dotsc$ & $1.54242401828716131644723995819\dotsc$  \\
$179$ & $3.86236132549903008112126130282\dotsc$ & $5.59074764196693719810304550344\dotsc$ & $1.60085064594072009293300914735\dotsc$  \\
$181$ & $5.14111848776848135810136664257\dotsc$ & $5.52401113238735460988935254057\dotsc$ & $1.65656567095010010041093792977\dotsc$  \\
$191$ & $4.69286990201422664003552434812\dotsc$ & $6.21621633683078754687889560801\dotsc$ & $1.69400806335478035992195123369\dotsc$  \\
$193$ & $5.16342219673915483320078262720\dotsc$ & $6.33516880970302226248749231989\dotsc$ & $1.72106839151430000218016220949\dotsc$  \\
$197$ & $7.55148715896640647886485129372\dotsc$ & $6.72431280547758930911931614898\dotsc$ & $1.58425224704856913591906318269\dotsc$  \\
$199$ & $6.47366513609320738699497459778\dotsc$ & $4.97867314026834059118807347477\dotsc$ & $1.52055512030192431037107983792\dotsc$  \\
$211$ & $7.73613578424586162532810587585\dotsc$ & $5.43928767077706865027592727891\dotsc$ & $1.58887689723521687477342354947\dotsc$  \\
$223$ & $7.81777971785991367471336734851\dotsc$ & $6.97640718267880419790301145060\dotsc$ & $1.57809439787964273689310796956\dotsc$  \\
$227$ & $8.08053156951296218697071193757\dotsc$ & $6.16478105833535800088839052312\dotsc$ & $1.61440476278289514090073256762\dotsc$  \\
$229$ & $7.16298632058099546745778115058\dotsc$ & $5.19368182825228459062582716349\dotsc$ & $1.64391627222705529854073112016\dotsc$  \\
$233$ & $3.11948354485127541303115295258\dotsc$ & $5.48268694035180653761326391137\dotsc$ & $1.56534808865669695863593307680\dotsc$  \\
$239$ & $3.99911017207833249512632297919\dotsc$ & $4.89826038220509731091188200357\dotsc$ & $1.83593237895342242137799671838\dotsc$  \\
$241$ & $6.03752521401034215065709250935\dotsc$ & $6.91099570349028181262249488655\dotsc$ & $1.74483502309356231328685290592\dotsc$  \\
$251$ & $5.04313708502347351042811119022\dotsc$ & $5.85522475367262429906377535883\dotsc$ & $1.60634233356394595761434310531\dotsc$  \\
$257$ & $8.16991391232741391670225155227\dotsc$ & $7.41413126491779482941571986652\dotsc$ & $1.52986363395322517571321794433\dotsc$  \\
$263$ & $7.30343624736815435414348077406\dotsc$ & $6.88761891078185993452639437420\dotsc$ & $1.61873689910065712561008039262\dotsc$  \\
$269$ & $6.26034831666577102735252755712\dotsc$ & $6.33572466741282346876839833227\dotsc$ & $1.58662353583078976012953348699\dotsc$  \\
$271$ & $5.97717804854803304223773905976\dotsc$ & $4.91607375378349595312704873315\dotsc$ & $1.51145118046000075647340279932\dotsc$  \\
$277$ & $4.59280817714077895164777081661\dotsc$ & $6.07306330239530923314413596279\dotsc$ & $1.72974155675277125427451583060 \dotsc$  \\
$281$ & $4.66496432366211457505220852623\dotsc$ & $4.99043740542558229612252801406\dotsc$ & $1.60536366070704717918242357661\dotsc$  \\
$283$ & $7.15028579741068251409225231188\dotsc$ & $7.04969230270522888347459792033\dotsc$ & $1.55609186296142373233316514603\dotsc$  \\
$293$ & $3.38438152121953978658468259238\dotsc$ & $5.38438152121953978658468259238\dotsc$ & $1.58515317244284064528356780036\dotsc$  \\
\hline
\end{tabular}
}
\caption{\label{table2}
Values of $\G_q$, $\G_q^+$ and $M_q$ for every odd prime up to $300$ with 
a precision of $30$ digits; computed with PARI/GP, v.~2.11.4
with  trivial summing over $a$. Total computation time: for $\G_q$, $\G_q^+$:
18 sec. 852 millisec.,  for $M_q$: 19 sec., 171 millisec. on the Dell
Optiplex machine mentioned before.}
\end{table}

\begin{table}[htp]
\begin{center}
\scalebox{0.7}{
\begin{tabular}{|r|r|r|r|r|r|} 
\hline
  $q$ &  $\G_{q}$\hskip3truecm\mbox{}  &$\G_{q}^+$ \hskip3truecm\mbox{}  &  time  &  time  & time  \\  
  &   && $T$ & $S$ & direct \\ \hline
$1009$ & $8.4421351518492992758606946727\dotsc$  & $6.2733540844322103172186250111\dotsc $& 5s. &   3s. & 14s.\\ 
$2003$ &  $5.7934213690793633280384982162\dotsc$  & $6.9935258611413978746616842142\dotsc$  & 10s. & 7s. & 39s.\\ 
$3001$ &  $8.6474651369683869388023453509\dotsc$  &  $8.6459700672984138998934976577\dotsc$ & 17s. & 11s. & 1m. 11s. \\ 
$4001$ &  $7.0034355462031439943568517684\dotsc$  &  $8.7805380094230735872867993849\dotsc$ & 24s. & 17s. & 1m. 49s. \\ 
$5003$ &  $5.5492930045816142277368795404\dotsc$  &  $7.2440224742791062634412330617\dotsc$ & 32s. & 23s. &  2m. 36s.\\ 
$6007$ &  $8.3116101219984838165629034403\dotsc$  & $9.8742666472425769486896123420\dotsc$  & 41s. & 30s. &  3m. 22s.\\ 
$7001$ &  $8.5052778761008771393168780384\dotsc$  &  $9.6833327734910786447084880544\dotsc$ & 52s. & 38s. &  4m. 07s.\\ 
$8009$ & $11.6868463915493575353450869960\dotsc$  &  $11.4431421556247084876087109206\dotsc$ & 1m. 03s. & 47s. &  5m. 00s.\\ 
$9001$ &  $10.1094784318383409358225035802\dotsc$  & $9.4868388831454962767492760006\dotsc$  & 1m. 15s. & 57s. &  5m.  56s. \\ 
$10007$ & $12.6646120045606923275389356783\dotsc$  & $11.0601624759024741933308283063\dotsc$  & 1m.  27s. & 1m. 07s.&  7m. 12s.\\
$20011$ & $10.7996803112999205186430402899\dotsc$  &  $10.5489807692170969459672226221\dotsc$ &  4m.  30s. &   3m. 43s.&  20m.  01s.\\
$30011$ & $10.3330799721240242255136062255\dotsc$  & $ 11.0127039500540893278498877674\dotsc$  & 9m. 19s.&  8m. 11s.& 37m.  28s.\\ 
\hline
\end{tabular} 
}
\caption{\label{table3}
A few other values of $\G_q$  and  $\G_q^+$ with 
a precision of $30$ digits; computed with PARI/GP, v.~2.11.4
with  trivial summing over $a$
[m = minutes, s = seconds]. Computation performed on the Dell
Optiplex machine mentioned before.}
\end{center}
\end{table}  %

   {\normalsize
\begin{table}[htp]
\begin{center}
\begin{tabular}{|r|r|r|r|} 
\hline 
 $q$\hskip0.5truecm\mbox{}   & $\G_{q}$  \hskip1truecm\mbox{}  & $\G_{q}^+$  \hskip1truecm\mbox{} &  time   \\\hline
$10007$ & $12.664612\dotsc$  & $11.060162\dotsc $ & 10ms. \\
$20011$ & $10.799680\dotsc$ & $10.548981\dotsc$ & 23ms. \\
$30011$ & $10.333080 \dotsc$ & $11.012704\dotsc $ & 15ms.    \\  
$40009$ & $13.146885\dotsc$ & $13.469520\dotsc$ & 25ms.  \\
$42611$ & $ 2.499688\dotsc$ & $8.367404\dotsc$ & 41ms.   \\  
$50021$ & $9.910507\dotsc$ & $11.063741\dotsc$ & 98ms.$^*$ \\ 
$60013$ & $12.810360\dotsc$ & $12.671109\dotsc$ & 36ms.   \\ 
$70001$ & $12.572765\dotsc$ & $13.428551\dotsc$ & 25ms. \\ 
$80021$ & $ 14.185633\dotsc$ &$ 11.617216\dotsc$ &  100ms.$^*$\\ 
$90001$ & $11.819424\dotsc$ & $9.601757\dotsc$ & 33ms. \\ 
$100003$ & $15.166074\dotsc$ & $14.765926\dotsc$ & 69ms. \\   
$305741$ & $1.650523\dotsc$ & $8.839799\dotsc$ & 198ms. \\ 
$1000003$ & $17.379970 \dotsc$ & $15.298449\dotsc$ &876ms. \\
$4178771$ & $ 0.922855\dotsc$ & $8.909168\dotsc$ &2s. 613ms. \\
$6766811$ & $1.604045\dotsc$ & $10.961044\dotsc$ &4s. 584ms. \\
$10000019$ & $17.087945\dotsc$ & $15.974742\dotsc$ &6s. 361ms. \\
$28227761$ & $2.361562\dotsc$ & $10.153369\dotsc$ &17s. 996ms. \\
$75743411$ & $2.469939\dotsc$ & $ 12.234097\dotsc$ &2m. 24s.  217ms. \\
  \hline 
\end{tabular}  
\caption{\label{table4}
A few other values of $\G_q$ and $\G_q^+$; computed with PARI/GP, v.~2.11.4 and {\tt fftw}, v.~3.3.8, with long double  precision.  The sum over $a$ was performed using the FFT algorithm on the Dell Optiplex machine mentioned before
[s = seconds, ms = milliseconds; precomputations of decimated in frequency $S$-values performed on the Optiplex; their computation  time is excluded from this table]. \\
$^*$: on the Intel Xeon machine due to a runtime  memory error on the Dell Optiplex.}
\end{center}
\end{table} 
 }

   {\normalsize
\begin{table}[htp]
\begin{center}
\begin{tabular}{|r|r|r|r|} 
\hline 
$q$ \hskip0.5truecm\mbox{}  &  $\G_{q}$  \hskip1truecm\mbox{} & $\G_{q}^+$  \hskip1truecm\mbox{} & time\\  \hline
$193894451$ & $0.662110\dotsc$ &$9.607705\dotsc$& 4m.  29s. \\ 
$212634221$& $1.435141\dotsc$ &$11.883540\dotsc$& 4m.  28s. \\
$251160191$& $ 1.912681\dotsc$ &$ 11.785574\dotsc$&2m. 53s. \\ 
$538906601$ & $1.474911\dotsc$ &$12.957235\dotsc$&11m.  56s.\\ 
$\pmb{964477901}$ & $\pmb{-0.182374\dotsc}$&$\pmb{10.402224}\dotsc$& \textbf{23m. 13s.}\\
$1139803271$ &$ 0.768538\dotsc$ &$8.313111\dotsc$ & 27m. 56s.  \\
$1217434451$ & $0.877596\dotsc$ &$12.946690\dotsc$ & 29m. 16s.  \\
$1806830951$ & $0.880396\dotsc$ &$11.973128\dotsc$& 47m. 48s. \\
$2488788101$ & $0.424880\dotsc$ &$ 12.248837\dotsc$& 103m. 08s. \\
$2830676081$ & $1.254528\dotsc$ &$12.438044\dotsc$& 89m. 59s. \\
$2918643191$ & $ 0.302793\dotsc$ & $12.573983\dotsc$& 87m. 49s. \\ 
$7079770931$ & $1.544698\dotsc$ & $14.301772\dotsc$& 742m. 09s. \\ 
$\pmb{9109334831}$ & $ \pmb{-0.248739\dotsc}$ &$\pmb{12.128187\dotsc}$& \textbf{311m.   28s.} \\
$\pmb{9854964401}$ & $ \pmb{-0.096465\dotsc}$ &$\pmb{12.807752\dotsc}$& \textbf{326m.   03s.} \\
 \hline 
\end{tabular} 
\caption{\label{table5}
A few other values of $\G_q$  and $\G_q^+$; computed with PARI/GP, v.~2.11.4 and {\tt fftw}, v.~3.3.8, with long double precision. 
Boldfaced results are the ones corresponding to known instances of $\G_q<0$.
The sum over $a$ was performed using the FFT algorithm on the Intel Xeon machine  or, for  $q=251160191$, $212634221$, $1139803271$, $7079770931$, $9109334831$, $9854964401$ on
the CAPRI infrastructure  mentioned before. 
%In particular the computational
%time for $\G_{9109334831}$ on CAPRI  was $311$ minutes and $28$ seconds;
%the computational time  for this case on the Xeon machine, as reported in this table, \textbf{1000m.   45s.} 
%was also  affected by  a runtime RAM swapping phenomenon.
[m = minutes, s = seconds; precomputations of decimated in frequency $S$-values performed on the cluster; 
their computation  time is excluded from this table].}
\end{center}
\end{table} 

{\normalsize
\begin{table}[htp]
\begin{tabular}{|r|r|}  
\hline
$n$  &  $\gamma_n$   \hskip4.5truecm\mbox{}  \\ \hline
$0$ & $0.5772156649015328606065120900824024310\dotsc$\\
$1$ & $-0.0728158454836767248605863758749013191\dotsc$\\
$2$ & $-0.0096903631928723184845303860352125293\dotsc$\\
$3$ &  $0.0020538344203033458661600465427533842\dotsc$\\
$4$ &  $0.0023253700654673000574681701775260680\dotsc$\\
$5$ &  $0.0007933238173010627017533348774444448\dotsc$\\
$6$ &  $-0.0002387693454301996098724218419080042\dotsc$\\
$7$ &  $-0.0005272895670577510460740975054788582\dotsc$\\
$8$ &  $-0.0003521233538030395096020521650012087\dotsc$\\
$9$ &  $-0.0000343947744180880481779146237982273\dotsc$\\
$10$ &  $0.0002053328149090647946837222892370653\dotsc$\\
$11$ &  $0.0002701844395439035266729020820679556\dotsc$\\
$12$ &  $0.0001672729121051401933535015433411834\dotsc$\\
$13$ &  $-0.0000274638066037601588600076036933551\dotsc$\\
$14$ &  $-0.0002092092620592999458371396973445849\dotsc$\\
$15$ &  $-0.0002834686553202414466429344749971269\dotsc$\\
$16$ &  $-0.0001996968583089697747077845632032403\dotsc$\\
$17$ & $0.0000262770371099183366994665976305101\dotsc$\\
$18$ &  $0.0003073684081492528265927547519486256\dotsc$\\
$19$ &  $0.0005036054530473556290555964377171600\dotsc$\\
$20$ &  $0.0004663435615115594494005948244335505\dotsc$\\
$21$ &  $0.0001044377697560001158107956743677204\dotsc$\\
$22$ &  $-0.0005415995822039977016551961731741055\dotsc$\\
$23$ &  $-0.0012439620904082457792997415995371658\dotsc$\\
$24$ &  $-0.0015885112789035615619061966115211158\dotsc$\\
$25$ &  $-0.0010745919527384888247242919873531730\dotsc$\\
$26$ &  $0.0006568035186371544315047730033562152\dotsc$\\
$27$ &  $0.0034778369136185382090073595742588115\dotsc$\\
$28$ &  $0.0064000685317006294581072282219458636\dotsc$\\
$29$ &  $0.0073711517704722391344124024235594021\dotsc$\\
$30$ &  $0.0035577288555731609479135377489084026\dotsc$\\ 

\hline
\end{tabular} 
\caption{\label{table6}
Computation of the generalised Euler constants $\gamma_n$, $0\le n\le 30$, with 
a precision of at least  40 digits; computed with PARI/GP, v.~2.11.4.}
\end{table}  
 
\newpage 
\begin{figure} [ht]
 \includegraphics[scale=0.94,angle=90]{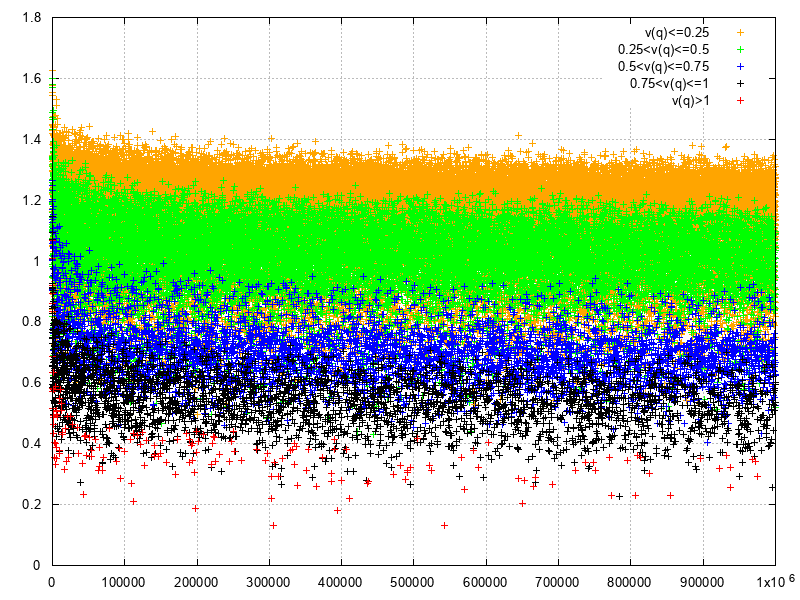}  
  \caption{The values of $\G_q/\log q$, $q$ prime, $3\le q\le 10^6$, plotted using
 GNUPLOT, v.5.2, patchlevel 8. 
 The  minimal value is $0.13067\dotsc$ and it is attained at $q=305741$; 
  the  maximal value is $1.62693\dotsc$ and it is attained at $q=19$. 
  Orange points satisfy $v(q)\le 0.25$; 
  green points satisfy $0.25 <v(q)\le 0.5$;    blue points satisfy $0.5 <v(q)\le 0.75$;
    black points satisfy $0.75 <v(q)\le 1$; red points satisfy $ v(q)> 1$; $v(q)$ is defined in \eqref{vq-def}. 
     }
 \label{fig1}
 \end{figure}  
\vfill\eject 
\begin{figure}[ht]
 \includegraphics[scale=0.94,angle=90]{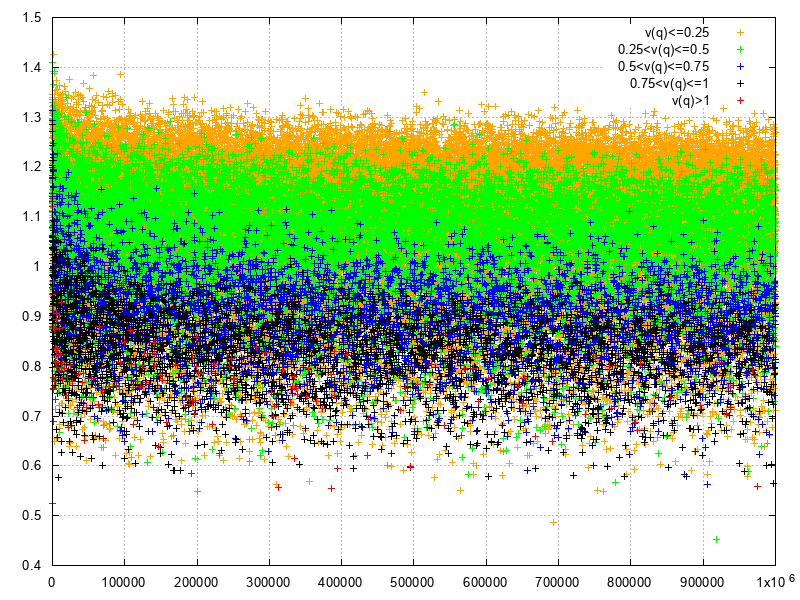}  
   \caption{The values of $\G_q^+/\log q$, $q$ prime, $3\le q\le 10^6$, plotted using
 GNUPLOT, v.5.2, patchlevel 8. 
 The  minimal value is $0.451468\dotsc$ and it is attained at $q= 918787$; 
 the  maximal value is $1.42626\dotsc$ and it is attained at $q=2053$.  
 Orange points satisfy $v(q)\le 0.25$; 
  green points satisfy $0.25 <v(q)\le 0.5$;    blue points satisfy $0.5 <v(q)\le 0.75$;
    black points satisfy $0.75 <v(q)\le 1$; red points satisfy $ v(q)> 1$; $v(q)$ is defined in \eqref{vq-def}.
    }  \label{fig2}
 \end{figure}
 
 \vfill\eject 
\begin{figure}[ht] 
 \includegraphics[scale=0.94,angle=90]{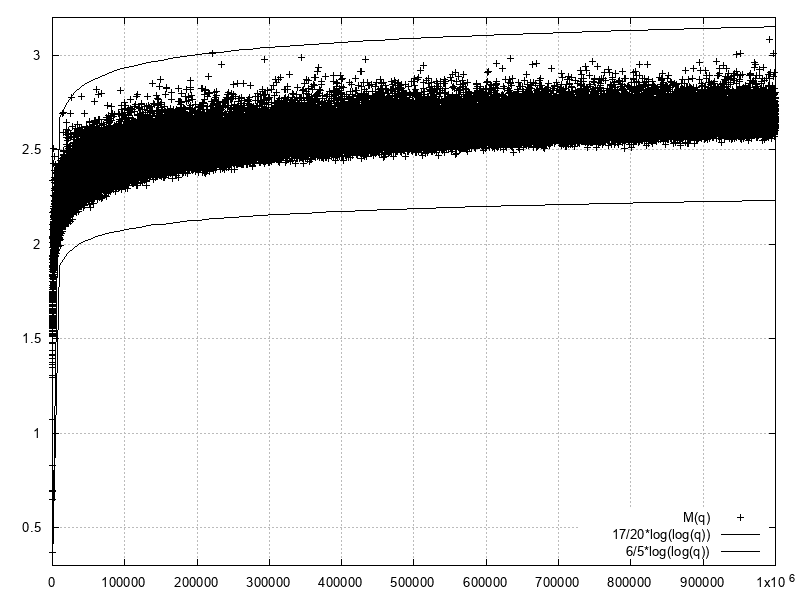}  
   \caption{The values of $M_q$, $q$ prime, $3\le q\le 10^6$, plotted using
 GNUPLOT, v.5.2, patchlevel 8. 
 The  minimal value is $0.3682816 \dotsc$ and it is attained at $q=3 $; 
 the  maximal value is $3.085536\dotsc$ and it is attained at $q= 991027$. 
 The lines  represent  the functions $c\cdot\log \log q$, with $c=17/20$, respectively  $c=6/5$. 
 }
  \label{fig3} 
 \end{figure}
  
   \vfill\eject 
\begin{figure}[ht] 
 \includegraphics[scale=0.94,angle=90]{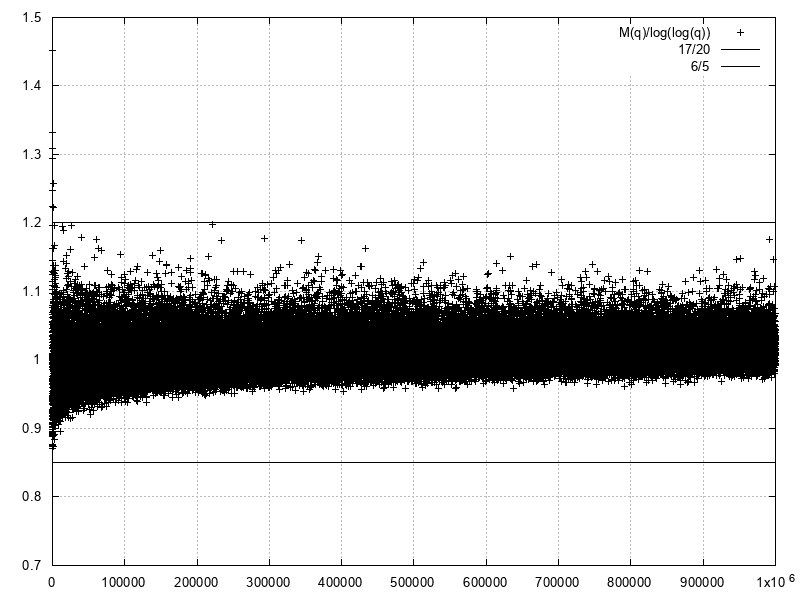}  
   \caption{The values of $M^\prime_q:=M_q/\log \log q$, $q$ prime, $3\le q\le 10^6$, plotted using
 GNUPLOT, v.5.2, patchlevel 8. 
 The  minimal value is $0.7392305\dotsc$ and it is attained at $q=13$; 
 the  maximal value is $3.9158971\dotsc$ and it is attained at $q= 3$ (not represented
 in the plot). 
 $M^\prime_q> 17/20$ for every $13< q\le 10^6$;
  $M^\prime_q <6/5 $ for every $1531< q\le 10^6$.
 The  lines  represent  the constant functions  $c=17/20$ and $c=6/5$. 
 }
  \label{fig4} 
 \end{figure}
  
\end{document}